# Outliers in dynamic factor models[*]

**Roberto Baragona**

*Dipartimento di Sociologia e Comunicazione, Sapienza Università di Roma, Via Salaria 113, I-00198 Roma, Italy*
*e-mail:* `roberto.baragona@uniroma1.it`

**Francesco Battaglia**

*Dipartimento di Statistica, Probabilità e Statistiche Applicate, Sapienza Università di Roma, Piazzale Aldo Moro 5, I-00100 Roma, Italy*
*e-mail:* `francesco.battaglia@uniroma1.it`

**Abstract:** Dynamic factor models have a wide range of applications in econometrics and applied economics. The basic motivation resides in their capability of reducing a large set of time series to only few indicators (factors). If the number of time series is large compared to the available number of observations then most information may be conveyed to the factors. This way low dimension models may be estimated for explaining and forecasting one or more time series of interest. It is desirable that outlier free time series be available for estimation. In practice, outlying observations are likely to arise at unknown dates due, for instance, to external unusual events or gross data entry errors. Several methods for outlier detection in time series are available. Most methods, however, apply to univariate time series while even methods designed for handling the multivariate framework do not include dynamic factor models explicitly. A method for discovering outliers occurrences in a dynamic factor model is introduced that is based on linear transforms of the observed data. Some strategies to separate outliers that add to the model and outliers within the common component are discussed. Applications to simulated and real data sets are presented to check the effectiveness of the proposed method.

## Contents

[*]Financial support by Sapienza University of Rome is gratefully acknowledged.

## 1. Introduction

Dynamic factor models have been introduced to explain and forecast time series of interest in the presence of a large set of explanatory time series. In practice, usefulness of dynamic factor models is apparent when the dimension is so large that vector autoregressive models are not able to handle the multivariate time series efficiently. Reduction of the available time series to few factors allows efficient and interpretable models to be estimated. Factor extraction has to be accomplished in such a way that only negligible or little amount of information be lost.

The study of the eigenvalues and eigenvectors of the parameter matrices was early suggested by (18) to produce a simplified version of an autoregressive model. A canonical transformation of a vector autoregressive model based on the simultaneous relationships between variables was introduced by (2) . The relationships between different time lags were considered by (9) and (19) in the frequency domain. The principal component analysis was extended in the frequency domain by (3). Identification of the number of factors in multivariate time series process was addressed by (14).

Factor models are strictly related to the diffusion indexes methodology (for instance 20). As pointed out by (6), when the dimension is large vector autoregressive (VAR) and vector autoregressive moving average (VARMA) models are difficult to estimate because the number of parameters grows with the number of time series quadratically. On the contrary, for dynamic factor models the growth is linear.

Usually mutually uncorrelated factors are assumed, whilst individual factor time series may be autocorrelated. The multivariate dynamic structure of the observed $N$-component time series $y_t$ may be modeled through the matrix factor $A$, that is

$$y_t = A x_t + \eta_t$$

where $x_t$ is a vector of $K$ independent time series and $\eta_t$ is the idiosyncratic disturbance. Each factor time series may follow a linear model, that is

$$x_{i,t} = \theta_i(B) \varepsilon_{i,t}$$

where $B$ is the back-shift operator and $\varepsilon_{i,t}$ is uncorrelated white noise. This leads to the dynamic factor model

$$y_t = A\theta(B)\varepsilon_t + \eta_t$$

where $\theta(B) = \text{diag}(\theta_1(B), \ldots, \theta_K(B))$. This model is a special case of

$$y_t = \psi(B)\varepsilon_t + \eta_t$$

as considered in (5), for instance. The equality

$$\psi_{i,j}(B) = a_{i,j}\theta_j(B)$$

reduces to the assumption that the impact of any shock $\varepsilon_{j,t}$ on the observed time series $y_{i,t}$ decays over time in similar way for any $i$. This assumption may also be justified on the ground of the asymptotic results by (6), p. 456.

Outliers in time series were introduced by (7) according to two different models, the additive outlier (or aberrant observation) and the innovation outlier (or aberrant innovation). This latter impacts the observed time series for some time span after the occurrence date, the former affects only one observation at the date of its occurrence. In spite of this, the additive outlier has serious effects on parameter estimates and forecasts, while the effects of the innovation outlier is often less serious. This motivates our choice for modeling outliers in the dynamic factor model as outlying observations of the additive type.

The plan of the paper is as follows. In Section 2 we introduce the outlier structure that we assume to be possibly present in a dynamic factor model. This structure and its implications will be examined in detail in Section 3. A method for checking the adequacy of the dynamic factor model to fit the data will be illustrated in Section 4. Methods for estimating the dynamic factor model parameters are discussed in Section 5 and the impact of outliers on the estimates will be examined in Section 6. In Section 7 a procedure for identifying outliers and estimating their size is presented and illustrated by an example. A simulation experiment for checking the effectiveness of the procedure in comparison with a multivariate model-based method (23) and a projection pursuit-based procedure (8) will be reported in Section 8. Our procedure is then applied to a set of real data, that is some quarterly economic data on asset prices, activity, wages, goods and commodity prices from the seven-country data set studied by (22). Results are reported in Section 9. Section 10 concludes. Proof of theorems are found in Appendix A.

## 2. The dynamic factor model with outliers

Let $y_t$ be an observed $N$-component vector time series and the temporal index $t = 1, \ldots, T$. We may even assume that the number of the time series components $N$ is greater than the number $T$ of dates when observations were made. We assume further that though $N$ may be very large the observed time series is actually explained by a much smaller number $K$ of unobservable time series $x_t = (x_{1t}, \ldots, x_{Kt})'$ and an idiosyncratic $N$-dimension disturbance $\eta_t$. Then the dynamic factor model with outliers may be written

$$y_t = Ax_t + \omega\Delta_t + \eta_t, \tag{2.1}$$

where $A$ is an $N \times K$ matrix of rank $K$, $K \ll N$. The outliers occurrence dates are modeled by the binary series $\{\Delta_t\}$ and by the $N \times 1$ vector $\omega$ which represents the outlier size.

Let us make about model (2.1) the following assumption:

1. $\{x_{1t}\}$, $\{x_{2t}\}$, ..., $\{x_{Kt}\}$ are mutually independent standardized random processes, i.e. $\mathbb{E}(x_{it}) = 0$, $\mathbb{E}(x_{it}^2) = 1$ for any $i$ and $t$, and $\mathbb{E}(x_{it} x_{js}) = 0$ if $i \neq j$.
2. The dynamics of the unobserved factor time series $x_t$ may be modeled as

$$x_{it} = \theta_{ii}^{(0)} \varepsilon_{it} + \theta_{ii}^{(1)} \varepsilon_{i,t-1} + \theta_{ii}^{(2)} \varepsilon_{i,t-2} + \ldots,$$

where $\theta_{ii}^{(0)} = 1$, $\sum_{j=0}^{\infty} (\theta_{ii}^{(j)})^2 < \infty$, $i = 1, \ldots, K$, and $\varepsilon_t = \{\varepsilon_{1t}, \ldots, \varepsilon_{Kt}\}$ are Gaussian white noises mutually independent at all leads and lags with diagonal variance-covariance matrix $\Sigma_\varepsilon$.
3. $\{\Delta_t\}$ is a deterministic scalar binary time series and $\omega$ is a non random $N \times 1$ vector. For an outlier occurring at time $t_0$, $\Delta_t = 1$ if $t = t_0$ and $\Delta_t = 0$ otherwise.
4. $\eta_t = \{\eta_{1t}, \ldots, \eta_{Nt}\}$ are Gaussian stationary time series both serially and mutually independent at all leads and lags with diagonal variance-covariance matrix $\Sigma_\eta$.
5. The vector time series $\varepsilon_t$ and $\eta_t$ are mutually independent at all leads and lags.

These assumptions are motivated by the idea that the dependence among the observed time series components is entirely explained by the factors. Therefore the idiosyncratic terms are also independent, since otherwise they would contribute to explain correlations between two observed components and should be put into the factor vector.

In general model (2.1) is not identified unless some assumptions are made about either the matrix $A$ or the vector time series $x_t$. In fact, let $C$ be any non-singular $K \times K$ matrix. Then in model (2.1) we have

$$A x_t = A C^{-1} C x_t.$$

By letting $A^* = AC^{-1}$ and $x_t^* = Cx_t$ model (2.1) could be written as well

$$y_t = A^* x_t^* + \omega \Delta_t + \eta_t.$$

No restriction is made on the matrix $A$ except that its rank is equal to $K$. Notice that Assumption 1 does not imply any loss of generality. In fact if $\Gamma_x(0) = \text{cov}(x_t, x_t')$ is not the identity matrix $I_K$ we could replace $x_t$ with the transformed data $\Gamma_x(0)^{-\frac{1}{2}} x_t$. As $\Gamma_x(0)$ is positive definite then a factorization $\Gamma_x(0) = \Gamma_x(0)^{\frac{1}{2}} (\Gamma_x(0)^{\frac{1}{2}})'$ exists for a non-singular matrix $\Gamma_x(0)^{\frac{1}{2}}$. The variance-covariance matrix of the transformed data turns out to be the identity matrix.

We prove the following theorem in Appendix A.1.

**Theorem 2.1.** *Model (2.1) under Assumptions 1 and 2 is identified up to factor sign changes.*

It has to be noticed that model (2.1) is uniquely determined by Assumptions 1—5 up to a permutation matrix and changing of sign. In fact, the order of the factors may be taken arbitrarily without affecting the model's structure. Moreover, Assumption 1 determines the factors sizes but each factor may be multiplied by $\pm 1$ without affecting its variance.

We may write the relationships that link the variance-covariance matrices of the observed data with those of the factors and of the idiosyncratic component

$$\Gamma_y(0) = A\Gamma_x(0)A' + \Sigma_\eta = AA' + \Sigma_\eta,$$

and

$$\Gamma_y(h) = A\Gamma_x(h)A', \ h \neq 0.$$

Let $\gamma_{i,j}^y(h)$, $i,j = 1,\ldots,N$, denote the entry in row $i$ and column $j$ of the matrix $\Gamma_y(h)$, and $\gamma_i^x(h)$, $i = 1,\ldots,K$, denote the diagonal elements of the matrix $\Gamma_x(h)$.

Note that Assumption 1 is used that implies $\Gamma_x(0) = I$ in the first equality, while the second equality shows that the matrices $\{\Gamma_y(h)\}$'s are symmetric because the $\{\Gamma_x(h)\}$'s are diagonal.

It is sometimes assumed that the columns of $A$ are orthogonal. This ensures the advantage that those columns are eigenvectors of all the covariance matrices of $\{y_t\}$ at any lag (see, e.g., 14). However, we feel that such assumption, together with Assumption 1, is somewhat unrealistic and it will not be formulated here.

## 3. Outliers in factor models

The estimation of outliers in Equation (2.1) is greatly simplified if a linear transform of the data exists that may highlight the impact of outlying observations. If parameters in Equation (2.1) are assumed known, then, by taking the projection matrix $Z = I - A(A'A)^{-1}A'$, the following lemma is easily proved.

**Lemma 1.** *Let the $N \times K$ matrix $A$ be defined as in Equation (2.1). Then a $N \times N$ matrix $Z$ exists such that $ZA = 0$ ( 0 is the $N \times K$ zero matrix).*

Lemma 1 has interesting implications concerned with the outliers estimation in model (2.1). The matrix $Z$ projects the vectors of $\mathbb{R}^N$ into the space orthogonal to the space $V_A$ spanned by the columns of $A$. Let $V_A^\perp$ denote this orthogonal space. By letting $V$ be the space of the vectors in $\mathbb{R}^N$ we have $V = V_A \oplus V_A^\perp$. Any vector in $V$ may be written as the sum of a vector in $V_A$ and a vector in $V_A^\perp$. Then three cases may occur

1. $Z\omega = 0$. In this case $\omega \in V_A$, that is there exist coefficients $c_1, c_2, \ldots, c_K$ such that

$$\omega = a_1 c_1 + a_2 c_2 + \ldots + a_K c_K,$$

where $a_1, a_2, \ldots, a_K$ denote the columns of $A$. We may write $\omega = Ac$, where $c = (c_1, c_2, \ldots, c_K)'$. Then Equation (2.1) becomes

$$y_t = A(x_t + c\Delta_t) + \eta_t.$$

The outliers are entirely within the factors, that is the observed $y_t$ are affected by outliers only through the factors.

2. $Z\omega \neq 0$ and $A(A'A)^{-1}A'\omega = 0$. This means that $\omega \in V_A^\perp$. The outliers impact the observed $y_t$ but the factors are actually outlier free.
3. $Z\omega \neq 0$ and $A(A'A)^{-1}A'\omega \neq 0$. Then $\omega$ may be written as the linear combination of a basis in $V$ obtained by assuming the columns of $A$ as a basis in $V_A$ and a basis $M = [m_1, m_2, \ldots, m_{N-K}]$ in $V_A^\perp$. We have

$$\omega = \sum_{i=1}^{K} a_i c_i + \sum_{j=1}^{N-K} m_j \mu_j = Ac + M\mu = Ac + \zeta \quad (3.1)$$

for some coefficients vectors $c$ and $\mu$. Model (2.1) becomes

$$y_t = A(x_t + c\Delta_t) + \zeta \Delta_t + \eta_t.$$

The observed time series $y_t$ is affected by an outlier of size $\zeta$ that adds to the whole structure and an additive outlier $c_r$ in each factor $x_{r,t}$.

Cases (1) and (3) may be treated by estimating $x_t^+ = x_t + c\Delta_t$ as if it were actually the model factors. Univariate search may be performed on the estimated $x_t^+$ factors to discover outliers dates and estimating their sizes.

We underline that in case (1), when $Z\omega = 0$, the dynamic model pattern is not affected by any perturbation. This latter is only transmitted by the model from factors to observed data. In that case the projection method we propose here is unable to identify the outliers, and they can only be discovered estimating the factors and employing univariate outlier search.

In cases (2) and (3) detection and estimation of outliers that impact the observed $y_t$ directly may be performed on the transformed model

$$Zy_t = Z\zeta\Delta_t + Z\eta_t$$

where $\zeta \in V_A^\perp$. Note that in case (2) we have $\omega = \zeta$ while in case (3) Equation (3.1) holds so that $\omega \neq \zeta$. In case (3) the outlier size $\omega$ has to be estimated partly in dynamic factor model and partly in the factors.

A similar development applies if the following dynamic model, as proposed by (5), is assumed:

$$y_t = \sum_{u=1}^{s} \psi_u \varepsilon_{t-u} + \eta_t + \omega \Delta_t \quad (3.2)$$

where the $\psi_u$'s are $N \times K$ matrices, $\varepsilon_t$ is a $K$-dimensional completely white noise with variances equal to 1, $sK < N$ and Assumptions 3, 4, 5 above hold. Let $V_\psi$ denote the space spanned by the columns of the matrices $\{\psi_u, u = 1, \ldots, s\}$ and $\mathbb{R}^N = V_\psi \oplus V_\psi^\perp$, and $Z$ the projection matrix onto $V_\psi^\perp$. We have

$$Zy_t = Z\eta_t + (Z\omega)\Delta_t .$$

In this case also $\omega$ may be written (but not necessarily in an unique way) as

$$\omega = \psi_1 c_1 + \psi_2 c_2 + \ldots + \psi_s c_s + \zeta$$

where $\zeta \in V_\psi^\perp$, and the model may be expressed as follows:

$$y_t = \sum_{u=1}^{s} \psi_u(\varepsilon_{t-u} + c_u \Delta_t) + \eta_t + \zeta \Delta_t$$

which decomposes the effect of an outlier into two parts, one that perturbs the dynamic structure of the model by altering the effect of the past values of the factors ($c_u$) and the other one simply superimposed to the observation ($\zeta$).

Usually the matrix $A$ is unknown and we may apply the preceding procedure only by computing an estimate $\hat{A}$. The presence of the outlying observations themselves makes the estimation difficult and often unreliable. Under some additional assumptions the following theorems allow an alternative procedure to be entertained which does not require estimating $A$. In what follows, all eigenvectors are normalized, i.e. they are taken with modulus equal to one.

**Theorem 3.1.** *Let $y_t$ satisfy model (2.1) with Assumptions 1—5 and suppose that for each $j = 1, 2, \ldots, K$ there exists a lag $h_j \neq 0$ such that $\gamma_j^x(h_j) \neq 0$. Then $z'A = 0$ if and only if $z$ is eigenvector associated with a zero eigenvalue of each $\Gamma_y(h), h \neq 0$.*

Proof is in Appendix A.2. Note that the assumptions of Theorem 3.1 are not satisfied if one of the factors is white noise.

**Theorem 3.2.** *Let $y_t$ satisfy model (3.2) and Assumptions 3—5 and suppose in addition that*
*(i) rank($\psi_1$) = $K$*
*(ii) There exists a lag $k$, $2 \leq k \leq s$ such that rank($\psi_k$) = $K$*
*then $z'\psi_u = 0, u = 1, 2, \ldots, s$ if and only if $z$ is eigenvector associated with a zero eigenvalue of each $\Gamma_y(h), h = 1, 2, \ldots, s - 1$.*

Proof is in Appendix A.3.

We note that Theorem 3.1 does not hold in general for $h = 0$ since in that case $\Gamma_y(0) = AA' + \Sigma_\eta$. Nevertheless, if the idiosyncratic disturbances are homoscedastic then the following theorems hold.

**Theorem 3.3.** *Let $z$ be any eigenvector associated to the smallest eigenvalue of $\Gamma_y(0)$. Let us assume, additionally, that $\Sigma_\eta = \sigma^2 I$. Then $z \in V_A^\perp$, that is $z'A = 0$. The converse is also true.*

Proof is in Appendix A.4.
A similar result holds for the dynamic model (3.2).

**Theorem 3.4.** *Let $y_t$ satisfy model (3.2) and suppose that $\Sigma_\eta = \sigma^2 I$. Then if $z'\psi_u = 0, u = 1, \ldots, s$, $z$ is eigenvector associated to the smallest eigenvalue of $\Gamma_y(0)$, and the converse is also true.*

Proof is in Appendix A.5.
The preceding theorems suggest a procedure to compute a projection of the multivariate time series that allows potential outliers to be readily detected.

If the hypothesis of homoscedasticity is assumed, we may compute an estimate $\hat{\Gamma}_y(0)$ (possibly a robust estimate) of the variance-covariance matrix $\Gamma_y(0)$. Then consider the eigenvectors associated to the smallest eigenvalue of $\hat{\Gamma}_y(0)$ (the smallest eigenvalue may have multiplicity greater than one). Let $z$ be any such eigenvector, then, according to Theorem 3.3 and 3.4, for the univariate time series $z'y_t$ we have

$$z'y_t \approx z'\eta_t + (z'\omega)\Delta_t.$$

Any such projection of the multivariate time series may be analyzed by means of univariate methods to detect potential outlying observations. Then evaluate the outlier's size by assuming the dates of occurrence of outliers from univariate analysis and using estimation methods in the multivariate framework.

If the homoscedastic hypothesis may not be assumed, the same result is obtained using Theorem 3.1 and 3.2, by taking $z$ equal to the eigenvector associated to a zero eigenvalue of a $\hat{\Gamma}_y(h)$ for $h > 0$.

## 4. Factor model adequacy

A crucial point is whether the simple factor model (2.1) together with our Assumptions fits the data adequately. Increasing the number of factors $K$ does not solve the problem because not all processes may be represented by equation (2.1) for arbitrary $K$ and under Assumptions 1—5, since their autocovariance matrices have to be symmetric as seen before. This suggests that a measure of adequacy of the factor model might be obtained by evaluating the differences between the elements $(i, j)$ and $(j, i)$ of $\Gamma_y(h)$, or the autocorrelation matrix. Let $r_{ij}(h) = \gamma_{ij}^y(h)\{\gamma_{ii}^y(0)\gamma_{jj}^y(0)\}^{-1/2}$, and denote by $\hat{r}_{ij}(h)$ the corresponding estimate. If $\Gamma_y(h)$ is symmetric, using classical results (see, e.g., 17) we obtain that the difference $\hat{r}_{ij}(h) - \hat{r}_{ji}(h)$ is asymptotically normal with mean zero and variance

$$\text{var}\{\hat{r}_{ij}(h) - \hat{r}_{ji}(h)\}$$
$$\simeq \frac{2}{T} \sum_{u=-\infty}^{\infty} [r_{ii}(u)r_{jj}(u) - r_{ij}(u)^2 - r_{ii}(u)r_{jj}(u-2h) + r_{ij}(u)r_{ij}(u-2h)]$$

and it depends both on the cross-correlation and autocorrelation functions in a complicated way; furthermore, such differences are correlated for different indexes $i$ and $j$. Therefore the differences $\hat{r}_{ij}(h) - \hat{r}_{ji}(h)$ cannot be used in any plausible way to test the hypothesis that $\Gamma_y(h)$ is symmetric. However, a possible solution is found turning to the frequency domain, in an analogue way as proposed by (9) when estimating parameters of factor models.

In the frequency domain the symmetry of the covariance matrices $\Gamma_y(h)$ for any $h$ is equivalent to a real spectral density matrix for any frequency. Let

$$F(\lambda) = \frac{1}{2\pi} \sum_{h=-\infty}^{\infty} \Gamma_y(h) e^{-i\lambda h}$$

denote the spectral density matrix of $y_t$. We prove the following theorem in Appendix A.6.

**Theorem 4.1.** *If $\Gamma_y(h)$ is symmetric for any h, then $F(\lambda)$ is real for any $\lambda$ and vice versa.*

We shall therefore test the hypothesis that $F(\lambda)$ is real for any $\lambda$.

Let us define the Fourier transforms as the $N \times 1$ complex vectors

$$d_T(\lambda) = \frac{1}{\sqrt{2\pi T}} \sum_{t=1}^{T} y_t e^{-i\lambda t}.$$

From Theorem 4.4.1 of (3) it follows that the real $2N \times 1$ vector $[\text{Re}\, d_T(\lambda)', \text{Im}\, d_T(\lambda)']'$ converges as $T \to \infty$ to a normal random vector with mean zero and variance covariance matrix:

$$\frac{1}{2} \begin{bmatrix} \text{Re}F(\lambda) & -\text{Im}F(\lambda) \\ \text{Im}F(\lambda) & \text{Re}F(\lambda) \end{bmatrix}.$$

If the spectral matrix is real, $\text{Im}\,F(\lambda) = 0$, thus $\text{Re}\,d_T(\lambda)$ and $\text{Im}\,d_T(\lambda)$ are (as $T \to \infty$) independently identically distributed normal vectors with zero means and variance covariance matrix $\frac{1}{2}F(\lambda)$. Therefore the hypothesis of real spectral density is equivalent to the independence of two normal vectors and may be tested by means of likelihood ratio. However, only one observation would be available for each fixed $\lambda$. To overcome such difficulty, and to test reality for any $\lambda$, we use a device similar to (9).

Let $\lambda_j = 2\pi j/T$ denote the Fourier frequencies, and suppose that $T$ is sufficiently large so that for a set of frequencies $\{\lambda_j, a < j \leq b\}$ we can assume $F(\lambda_j) \simeq F$. Also, let $J = b - a$ and

$$X_j^R = \text{Re}\, d(\lambda_{a+j}) \quad X_j^I = \text{Im}\, d(\lambda_{a+j}) , \quad j = 1, \ldots, J$$

where we have dropped for convenience the dependence on $T$. For $T$ large we may assume that

$$\begin{pmatrix} X_j^R \\ X_j^I \end{pmatrix} \sim N\left[\begin{pmatrix} 0 \\ 0 \end{pmatrix}, \frac{1}{2}\begin{pmatrix} \text{Re}F & -\text{Im}F \\ \text{Im}F & \text{Re}F \end{pmatrix}\right]$$

while under the null hypothesis $H_0$: $F$ real, $X_j^R$ and $X_j^I$ are independently identically distributed normal vectors with zero means and variance covariance matrix $\frac{1}{2}F$. Define the variance estimates:

$$S_R = \frac{1}{J}\left(\sum_{j=1}^{J} X_j^R(X_j^R)' + X_j^I(X_j^I)'\right) , \quad S_I = \frac{1}{J}\left(\sum_{j=1}^{J} X_j^I(X_j^R)' - X_j^R(X_j^I)'\right).$$

The following theorem provides the likelihood ratio test.

**Theorem 4.2.** *The likelihood ratio test statistic for the null hypothesis $H_0$: $F$ real is given by*

$$\left|I + (S_R^{-1}S_I)^2\right|$$

*and its distribution under $H_0$ is equal to that of the statistic $U_{N,N,J-N-1}$ of (1, chap. 9).*

Proof is in Appendix A.7. Some approximations are discussed in (1), which for our statistic imply approximating $-m\log U$ ($m = J - N - \frac{3}{2}$) by a chi–square variable with $N^2$ degrees of freedom. The test of Theorem 4.2 may be employed repeatedly on non-overlapping frequency intervals, and the usual caveats for multiple testing apply (see, e. g., 11, chap. 5.4).

## 5. Estimation problems

Outliers identification is only concerned with the detection of time points where outlying observations occur. When this task is performed by examining a univariate projection series, as is suggested in this paper, little may be said about the multivariate outliers size. The model parameters matrices, either $A$ in the model's formulation (2.1) or $\{\psi_1, \psi_2, \ldots\}$ in the model's formulation (3.2), have to be estimated from available data $\{y_t, t = 1, \ldots, T\}$, along with the common components and the variance covariance matrix $\Sigma_\eta$ of the idiosyncratic component. This way outliers size may be estimated while the estimated model is available for studying the relations among the observed time series or for forecasting purpose.

We shall distinguish in estimation procedures whether the idiosyncratic covariance matrix is constrained to the relationship $\Sigma_\eta = \sigma^2 I$, or to be a diagonal matrix, or no special constraints are imposed on its entries. Also, we consider here the unperturbed case of absence of outliers. The distortion induced by the presence of an outlier will be considered in the next section.

The log-likelihood of $\{y_1, y_2, \ldots, y_T\}$ according to model (2.1) under the assumption $\Sigma_\eta = \sigma^2 I$ is

$$L = -\frac{NT}{2}\log 2\pi - NT\log\sigma^2 - \frac{1}{2\sigma^2}\sum_{t=1}^{T}(y_t - Ax_t)'(y_t - Ax_t).$$

Let us define the $N \times T$ matrix

$$Y = [y_1, y_2, \ldots, y_T],$$

where the $y_t$'s are $N \times 1$ arrays, and the $K \times T$ matrix

$$X = [x_1, x_2, \ldots, x_T],$$

where the $x_t$'s are $K \times 1$ arrays. Then the sum of squares in the log-likelihood may be written

$$\operatorname{tr}\{(Y - AX)'(Y - AX)\},$$

and its minimization is equivalent to maximizing the likelihood.

Model (2.1) is considered by (20) and (21) who assume normality, $\Sigma_\eta = \sigma^2 I$, and treat the factors $\{x_t\}$ as deterministic components. They show that, on maximizing the likelihood with respect to both $\{x_t\}$ and $A$, the maximum likelihood estimate of $A$ is given by the matrix formed by the $K$ eigenvectors associated to the $K$ largest eigenvalues of $\hat{\Gamma}_y(0)$. They assume $A'A = I$; if we want to dispense with such hypothesis, we may use instead the fact that $\Gamma_x(0) = I$ (actually $XX' = I$ is assumed for simplicity). This leads to the following different estimate.

**Theorem 5.1.** *Let $\{y_t\}$ satisfy model (2.1) with Assumption 4 with $\Sigma_\eta = \sigma^2 I$, and suppose that $\{x_t\}$ are constants and $XX' = I$. Let $r$ be the rank of $Y$ and $K$ be a known pre-specified integer, so that $0 < K \leq r \leq \min(N,T)$, and let $\lambda_1 \geq \lambda_2 \geq \ldots \geq \lambda_T$ be the eigenvalues of $Y'Y$ with associated orthogonal eigenvectors $u_1, u_2, \ldots, u_T$ in $\mathbb{R}^T$. Then the problem*

$$\min\{\operatorname{tr}((Y - AX)(Y - AX)')\} \quad \text{subject} \quad \text{to} \quad XX' = I_K$$

*is solved by*
$$\hat{X} = [u_1, u_2, \ldots, u_K]'$$

*and*
$$\hat{A} = Y\hat{X}'.$$

*Further, the formula for $\hat{A}$ reduces to*
$$\hat{A} = W_K \Lambda_K^{1/2},$$

*where $W_K = [w_1, \ldots, w_K]$ ($w_1$, $w_2$, …, $w_K$ are the eigenvectors associated to the K largest eigenvalues of $YY'$), and $\Lambda_K^{1/2} = \text{diag}(\sqrt{\lambda_1}, \ldots, \sqrt{\lambda_K})$.*

Proof is in Appendix A.8. The estimate of the matrix $A$ given by Theorem 5.1 is consistent as the estimate $YY'/T$ is known to be consistent and its eigenvalues and eigenvectors (which are continuous functions of the elements of the matrix $YY'/T$) are consistent as well. If the observed data are not standardized, and if the eigenvalues are all distinct and the true variance covariance matrix is definite positive, then it may be shown that the eigenvalues are asymptotically independently normally distributed. The difference between estimated eigenvalues and actual ones is of order $1/\sqrt{T}$ in probability. The estimates of the eigenvectors are asymptotically normally distributed but they are not independent (see, e.g., 16, p. 290). Rate of convergence to actual eigenvectors is of order $1/\sqrt{T}$ in probability.

If we assume that the factors $\{x_t\}$ are random processes, the method of linear factor analysis may be employed. To overcome the problem that the factors are autocorrelated, (9) has introduced a frequency domain extension of the factor analysis which may be directly applied to model (2.1) assuming that $\Sigma_\eta$ is diagonal but not necessarily homoscedastic.

An alternative estimation method is using a Kalman filter in a state space formulation of the model, where (2.1) is considered as a measurement equation and $\{x_t\}$ is the state vector. In that case, a transition equation has to be specified for the factors $x_t$ which may be convenient if we assume that the process $\{x_t\}$ is easily modeled in state space form (if, for example, it is assumed a low order autoregression).

Finally, an estimation method which does not rely on any assumption on $\Sigma_\eta$ may be obtained using a technique of temporal blind source separation, for instance the temporal decorrelation source separation method (25) which uses an algorithm for approximate simultaneous diagonalization of several covariance matrices. Under model (2.1) we have
$$\Gamma_y(h) = A\Gamma_x(h)A' \qquad h \neq 0,$$
and, taking the generalized inverse $B = (A'A)^{-1}A'$,
$$\Gamma_x(h) = B\Gamma_y(h)B'.$$

Since $\Gamma_x(h)$ is diagonal for any $h \neq 0$, we may determine $B$ in such a way that the off diagonal elements of $B\Gamma_y(h)B'$ are as small as possible. Formally, the matrix $B$ is obtained by the following approximation problem
$$\min_B \sum_{h=1}^{H} \sum_{i \neq j} (B\Gamma_y(h)B')_{ij}^2.$$

Let $\hat{B}$ denote the solution, then

$$\hat{A} = \hat{B}'(\hat{B}\hat{B}')^{-1}.$$

A maximum lag $H$ has to be chosen. This may be selected by estimating the covariance matrices of the data $\hat{\Gamma}_y(.)$ and taking $H$ the minimum lag such that all entries of $\hat{\Gamma}_y(h)$, for $h > H$, are not significantly different from zero. Moreover, in order to recover the matrix $A$ is necessary that $(\hat{B}\hat{B}')^{-1}$ exists, therefore the solution matrix $\hat{B}$ should have full rank. Though it is easily seen that under model (2.1) this method is consistent, its sample properties appear very hard to be devised.

If model (3.2) seems more suitable to describe the data, the estimation methods proposed by (5) may be applied.

## 6. Bias induced by the outliers on the estimates

We turn now to consider the bias induced on the estimates of time-domain and frequency-domain indexes by the presence of an outlier.

Suppose that the observed time series $\{y_t, t = 1, \ldots, T\}$ contains an outlier at time $t_0$ and size measured by the vector $\omega$. Let $z_t$ denote the unperturbed data,

$$z_t = \begin{cases} y_t & \text{if } t \neq t_0 \\ y_{t_0} - \omega & \text{if } t = t_0 \end{cases}.$$

We compare the autocovariance or spectral density estimates computed on the actually observed series with the corresponding unbiased estimates that would have been obtained from the unperturbed time series $z_t$.

Note first that for the average $\bar{y}$ we have

$$\bar{y} = \frac{1}{T}\sum y_t = \frac{1}{T}\sum z_t + \frac{1}{T}\omega = \bar{z} + \frac{1}{T}\omega,$$

therefore

$$y_t - \bar{y} = z_t - \bar{z} - \frac{1}{T}\omega + \omega\Delta_t,$$

where $\Delta_t = 1$ if $t = t_0$ and zero otherwise. If we denote

$$\tilde{\Gamma}(h) = \frac{1}{T}\sum_{t=1}^{T-h}(z_t - \bar{z})(z_{t+h} - \bar{z})'$$

the unperturbed autocovariance estimates (computed on the unobserved time series $z_t$), the actual estimates may be written

$$\hat{\Gamma}(h) = \frac{1}{T}\sum_{t=1}^{T-h}(y_t - \bar{y})(y_{t+h} - \bar{y})' = \frac{1}{T}\sum_{t=1}^{T-h}(z_t - \bar{z} - \frac{1}{T}\omega + \omega\Delta_t)(z_{t+h} - \bar{z} - \frac{1}{T}\omega + \omega\Delta_{t+h})',$$

and a simple calculation gives

$$\hat{\Gamma}(0) = \tilde{\Gamma}(0) + \frac{\omega\omega'}{T} + \frac{1}{T}\{(z_{t_0} - \bar{z})\omega' + \omega(z_{t_0} - \bar{z})'\} - \frac{\omega\omega'}{T^2}$$

and

$$\hat{\Gamma}(h) = \tilde{\Gamma}(h) + \frac{1}{T}\{(z_{t_0-h} - \bar{z})\omega' + \omega(z_{t_0+h} - \bar{z})'\}$$
$$- \frac{1}{T^2}\sum_{t=1}^{T-h}\{(z_t - \bar{z})\omega' + \omega(z_{t+h} - \bar{z})'\} - \frac{T+h}{T^3}\omega\omega'. \quad (6.1)$$

It follows that the difference between unperturbed and actual estimates is of order $1/T$ in probability. We prove in Appendix A.9 the following theorem which states a more precise result.

**Theorem 6.1.** *Let $\{y_t, t = 1,\ldots,T\}$ be part of a realization of a second order stationary process with finite covariance matrices and suppose that at time $t_0$ an outlier equal to $\omega$ is added. Then if the $\hat{\gamma}_{ij}(h)$'s denote the covariance estimates and the $\tilde{\gamma}_{ij}(h)$'s denote the correspondent estimates computed on the unperturbed time series,*

$$\mathbb{E}\{T(\hat{\gamma}_{rs}(0) - \tilde{\gamma}_{rs}(0))\} = \omega_r\omega_s + O(T^{-1})$$

$$\text{var}\{T(\hat{\gamma}_{rs}(0) - \tilde{\gamma}_{rs}(0))\} = \omega_s^2\gamma_{rr}(0) + \omega_r^2\gamma_{ss}(0) + 2\omega_s\omega_r\gamma_{rs}(0) + O(T^{-1}),$$

*and, for $h \neq 0$,*

$$\mathbb{E}\{T(\hat{\gamma}_{rs}(h) - \tilde{\gamma}_{rs}(h))\} = O(T^{-1})$$

$$\text{var}\{T(\hat{\gamma}_{rs}(h) - \tilde{\gamma}_{rs}(h))\} = \omega_s^2\gamma_{rr}(0) + \omega_r^2\gamma_{ss}(0) + 2\omega_s\omega_r\gamma_{rs}(2h) + O(T^{-1}).$$

Results in the frequency domain may also be obtained, with the usual asymptotic approximations. On denoting by

$$\tilde{d}(\lambda) = \frac{1}{\sqrt{2\pi T}}\sum_{t=1}^{T}(z_t - \bar{z})\exp(-i\lambda t)$$

the Fourier transforms of the unperturbed data, we may write

$$d(\lambda) = \frac{1}{\sqrt{2\pi T}}\sum_{t=1}^{T}(y_t - \bar{y})\exp(-i\lambda t) = \frac{1}{\sqrt{2\pi T}}\sum_{t=1}^{T}(z_t - \bar{z} + \omega\Delta_t - \frac{\omega}{T})\exp(-i\lambda t)$$

$$= \tilde{d}(\lambda) + \frac{1}{\sqrt{2\pi T}}\omega\exp(-i\lambda t_0) - \frac{1}{\sqrt{2\pi T}}\frac{\omega}{T}\sum_{t=1}^{T}\exp(-i\lambda t).$$

If we consider only the Fourier frequencies $\lambda_j = \frac{2\pi}{T}j$, the third term disappears and

$$d(\lambda_j) = \tilde{d}(\lambda_j) + \frac{1}{\sqrt{2\pi T}}\omega\cos(\lambda_j t_0) - i\frac{1}{\sqrt{2\pi T}}\omega\sin(\lambda_j t_0),$$

and simple calculations show that for the actual and unperturbed periodograms

$$I(\lambda) = d(\lambda)\overline{d(\lambda)}' \quad , \quad \tilde{I}(\lambda) = \tilde{d}(\lambda)\overline{\tilde{d}(\lambda)}'$$

we have

$$I(\lambda_j) = \tilde{I}(\lambda_j) + \frac{\omega\omega'}{2\pi T} + A(\lambda_j) + iB(\lambda_j),$$

where

$$A(\lambda) = \frac{1}{2\pi T} \sum_{t=1}^{T} \cos\lambda(t_0 - t)\{\omega(z_t - \bar{z})' + (z_t - \bar{z})\omega'\},$$

and

$$B(\lambda) = \frac{1}{2\pi T} \sum_{t=1}^{T} \sin\lambda(t_0 - t)\{(z_t - \bar{z})\omega' - \omega(z_t - \bar{z})'\}.$$

We prove in Appendix A.10 the following theorem which gives more specific results as regards each entry of the matrix of differences between periodograms.

**Theorem 6.2.** *Let $\{y_t, t = 1,\ldots,T\}$ be part of a realization of a second order stationary process with spectral density matrix $F(\lambda)$ and suppose that at time $t_0$ an outlier equal to $\omega$ is added to the time series. Then if $I(\lambda)$ denotes the periodogram matrix and $\tilde{I}(\lambda)$ denotes the periodogram matrix of the correspondent unperturbed series, for each $1 \leq r \leq s \leq N$,*

$$\mathbb{E}T\operatorname{Re}\{I_{rs}(\lambda_j) - \tilde{I}_{rs}(\lambda_j)\} = \frac{\omega_r \omega_s}{2\pi}$$

*and*

$$\mathbb{E}T\operatorname{Im}\{I_{rs}(\lambda_j) - \tilde{I}_{rs}(\lambda_j)\} = 0.$$

*Let $s(T)$ be a sequence such that $\lambda(T) = 2\pi s(T)$ tends to $\lambda$ as $T \to \infty$. Then, as $T \to \infty$,*

$$\operatorname{var}\sqrt{(T)}\operatorname{Re}\{I_{rs}(\lambda(T)) - \tilde{I}_{rs}(\lambda(T))\} \to \frac{1}{4\pi}\{\omega_r^2 f_{ss}(\lambda) + \omega_s^2 f_{rr}(\lambda) + \omega_r \omega_s (f_{rs}(\lambda) + f_{rs}(-\lambda))\}$$

*and*

$$\operatorname{var}\sqrt{(T)}\operatorname{Im}\{I_{rs}(\lambda(T)) - \tilde{I}_{rs}(\lambda(T))\} \to \frac{1}{4\pi}\{\omega_r^2 f_{ss}(\lambda) + \omega_s^2 f_{rr}(\lambda) - \omega_r \omega_s (f_{rs}(\lambda) + f_{rs}(-\lambda))\}.$$

The preceding results enable also to evaluate the bias induced by an outlier on spectral estimates. Let

$$\hat{F}(\lambda) = \frac{2\pi}{T} \sum_{j=-T/2}^{T/2} I(\lambda_j) w_T(\lambda - \lambda_j)$$

be a consistent spectral estimate where $\lambda_j = \frac{2\pi}{T}j$ and $w_T(.)$ is a spectral window (see, e.g., 17), and denote by $\tilde{F}(\lambda)$ the same form computed on $\tilde{I}(\lambda)$. We assume that for each $T$ a truncation point $M = M(T)$ has been selected such that

$$\frac{2\pi}{T} \sum_j w_T(\lambda_j)^2 = c_0 M$$

and, as $T \to \infty$, $M \to \infty$ but $M/T \to 0$.

**Theorem 6.3.** *Under our assumptions as $T \to \infty$*

$$\mathbb{E}T\{\hat{F}_{rs}(\lambda) - \tilde{F}_{rs}(\lambda)\} = \frac{\omega_r \omega_s}{2\pi},$$

$$\text{var}\{\frac{T}{\sqrt{M}}\text{Re}(\hat{F}_{rs}(\lambda) - \tilde{F}_{rs}(\lambda) - \frac{\omega_r \omega_s}{2\pi})\}$$
$$\to \frac{c_0}{2}\{\omega_r^2 f_{ss}(\lambda) + \omega_s^2 f_{rr}(\lambda) + \omega_r \omega_s (f_{rs}(\lambda) + f_{rs}(-\lambda))\} \quad (6.2)$$

$$\text{var}\{\frac{T}{\sqrt{M}}\text{Im}(\hat{F}_{rs}(\lambda) - \tilde{F}_{rs}(\lambda))\} \to \frac{c_0}{2}\{\omega_r^2 f_{ss}(\lambda) + \omega_s^2 f_{rr}(\lambda) - \omega_r \omega_s (f_{rs}(\lambda) + f_{rs}(-\lambda))\}.$$

Proof is in Appendix A.11.

The preceding discussion offers some arguments for choosing among estimation methods. The methods based on the variance-covariance matrix or on the spectral density of the observed data have the drawback that the estimates may be biased by the presence of the outlier. The influence of the outlier on the estimates is of order $1/T$ and this would suggest that for large samples the two methods may perform equally well. Spectral methods, however, seem to constitute the favorite device to be used for checking dynamic factor model adequacy. The methods based on temporal decorrelation on the average are not affected by the presence of the outlier. This would suggest that these methods are likely to yield more reliable results. However a treatment of the sampling properties is not available and model adequacy has to be checked in any case to ensure that temporal decorrelation may be applied properly. Summing up, it seems that no estimation method may be declared the best one, and trying different methods seems advisable.

## 7. Outlier identification and size estimation procedure

Let $Y$ be the observed data arranged as a matrix of $N$ rows and $T$ columns. The entry $Y_{it}$ stands for the observation at time $t$ of the $i$-th time series, $t = 1, \ldots, T$ and $i = 1, \ldots, N$. We assume that a dynamic factor model may be tentatively fitted to the data with unknown number of factors and possibly outlying observations of unknown size at unknown dates. The idiosyncratic component is unknown as well. Given Assumptions 1—4 and the observed data set, model (2.1) is to be identified and estimated.

A procedure for estimating the unknown components of model (2.1) and performing outlier identification and estimation may be described as follows.

1. The optimal direction for locating the occurrence of outlying observations requires the computation of the variance-covariance matrix of the data $\hat{\Gamma}_y(0)$ and of its smallest eigenvalue $\lambda$, i.e.

$$\hat{\Gamma}_y(0) = \frac{1}{T}(Y - \bar{Y})(Y - \bar{Y})'$$

   where $\bar{Y} = \bar{y}1'$. The $N \times 1$ array $\bar{y}$ is the components average computed over time and 1 is the all-one $T \times 1$ vector. Then the eigenvalue $\lambda$ and the associated $N \times 1$ eigenvector $z$ may be computed. The linear transform $w = z'Y$ yields the $1 \times T$ time series that may be searched for the outlier occurrence dates. Let $\bar{w}$ and $\sigma_w$ be the mean and standard error of $w$. Then the outlier is identified at time $t$ if $|w_t|$ is the largest value such that $|w_t - \bar{w}| > \alpha \sigma_w$ where $\alpha$ is a suitable constant.

Since this procedure is based on Theorem 3.3, it strictly holds only in the homoscedastic case, and is only approximately appropriate when such hypothesis does not hold. An alternative would be computing the eigenvector associated to a zero eigenvalue of one of the matrices $\hat{\Gamma}_y(h)$ for $h > 0$, for example $\hat{\Gamma}_y(1)$, according to Theorem 3.1 and 3.2.

2. If an outlier is detected, then the multivariate time series may be adjusted by interpolating (by multivariate linear interpolator) or forecasting (by vector autoregressive (VAR) model) its value at the outlier date. Another strategy may consist in assuming that there is a missing value at the time of the outlier possible occurrence. Anyway, the potential outlier is either replaced by its conditional mean or a missing value is assumed, and the outlier free estimates of variance-covariance and lagged covariance matrices may be computed up to a pre-specified maximum lag $M$. Robust estimation methods may constitute a valuable choice. Either outlier free or robust estimates of the spectral density matrices have to be computed at the Fourier frequencies $\lambda_j$, where $\lambda_j = 2\pi j/T$ and $j = -T/2 + 1, -T/2 + 2, \ldots, T/2$.

3. Checking conformity of the observed data to the dynamic factor model may be done by using the spectral density estimates. For instance, let the frequency interval $0 < \lambda \leq \pi$ be divided in 4 sub-intervals and assume for simplicity that $T$ is integer multiple of 4. More than four sub-intervals may obviously be used, also according to the number of the data. If there is no reason for privileging any frequency components, equally wide sub-intervals may be selected. Symmetry considerations allow us to consider only the interval $(0, \pi)$ instead the whole interval $(-\pi, \pi)$. Then in each sub-interval $(\lambda_{a+1}, \lambda_b)$ there are $J = T/4$ Fourier frequencies, i.e. $\{\lambda_j, j = a+1, a+2, \ldots, a+J = b\}$. In the $\ell$-th sub-interval $a = (\ell-1)J$ and $b = \ell J$. The likelihood ratio test statistic $U$ is provided by Theorem 4.2. The null hypothesis $H_0$ that the covariance matrices are symmetric, i.e. the dynamic factor model may not be rejected, has to be checked by using the approximate statistic $-m \log U$ ($m = J - N - \frac{3}{2}$) which, under $H_0$, is distributed as a chi-square with $N^2$ degrees of freedom.

4. Once the model has been found appropriate, the number of factors $K$ has to be estimated. A simple device is based on the eigenvalues of the variance-covariance matrix $\tilde{\Gamma}_y(0)$ that is $\hat{\Gamma}_y(0)$ corrected for potential outliers. Let

$$v = \{v_1, v_2, \ldots, v_N\}$$

be the eigenvalues of $\tilde{\Gamma}_y(0)$ arranged in non increasing order and consider the cumulated sums

$$V_1 = v_1, \quad V_2 = V_1 + v_2, \quad \ldots, \quad V_N = V_{N-1} + v_N.$$

We assume $K$ as the number of factors if it is the smallest integer such that $V_K/V_N > 1 - \alpha$, where $\alpha$ is a small real number, 0.05 for instance. That is, the number of factors is chosen so as the cumulated normalized eigenvalues of $\tilde{\Gamma}_y(0)$ exceed a fraction $1 - \alpha$ that is judged large enough. We may want to take into account that the smallest eigenvalues may be greater than zero. So a better approximation to the correct number of factors may be obtained by assuming $K$ as the smallest integer such that $\{V_K + (N-K)\lambda_{\min}\}/V_N > 1-\alpha$. More information may be used for choosing $K$ along these same guidelines by using the eigenvalues of the (corrected) spectral density matrix at frequency zero which is essentially the sum of the covariance matrices for lags $-M, -M+1, \ldots, 0, 1, \ldots, M$. The covariance matrices may be used as well separately, and so the spectral matrices at non zero frequencies. However, different matrices may lead to different estimates of $K$ though the same value is to be expected in most cases.

5. Estimation of matrix $A$ and factors $X$ allows the dynamic factor model to be specified completely. In addition, both $A$ and $X$ are needed to estimate the outlier size and to distinguish if the factors, the model or both are affected by the outlying observations. Several methods are available, for example we may list the following ones.

   (a) Temporal decorrelation, i.e. a matrix $\hat{A}$ not necessarily squared nor orthogonal may be computed by approximate joint diagonalization (24) of matrices $\hat{\Gamma}_y(h)$, $h = 1, \ldots, K$. Let the $N \times K$ matrix $\hat{A}$ be such matrix, that is

   $$\hat{\Gamma}_y(h) = \hat{A}\hat{D}_h\hat{A}', \qquad h = 1, \ldots, K,$$

   where the $\hat{D}_h$'s are diagonal $K \times K$ matrices. Then we may let

   $$\hat{B} = (\hat{A}'\hat{A})^{-1}\hat{A}'$$

   and

   $$\hat{X} = \hat{B}Y.$$

   It follows
   $$\hat{B}\hat{\Gamma}_y(h)\hat{B}' = \hat{\Gamma}_x(h) = \hat{D}_h \qquad h = 1, \ldots, K$$

   for
   $$\hat{\Gamma}_x(h) = \frac{1}{T}\hat{X}\hat{X}' = \frac{1}{T}\hat{B}YY'\hat{B}' = \hat{B}\hat{\Gamma}_y(h)\hat{B}' = \hat{D}_h.$$

   (b) Assuming, for the purpose of estimating $A$, that the factors are not random, as suggested by (21; 22), the method outlined in Theorem 5.1 may be used, which only requires the $K$ eigenvectors associated to the $K$ largest eigenvalues of the matrix $YY'$ be computed. Assuming $W_K$ the $N \times K$ matrix whose columns are the column eigenvectors and $\Lambda_K$ the diagonal matrix with the largest eigenvalues of $YY'$ on its diagonal, the estimate of $A$ is given by the simple formula $\hat{A} = W_K\Lambda_K^{1/2}$.

   (c) Methods developed in the factor analysis framework may be used. Consider again the log-likelihood of the dynamic factor model

   $$L = -\frac{NT}{2}\log 2\pi - \frac{T}{2}\log|\Sigma_\eta| - \frac{1}{2}\text{tr}\{(Y-AX)'\Sigma_\eta^{-1}(Y-AX)\}.$$

   Let us assume at this stage that the matrix $\Sigma_\eta$ is known. Maximizing the likelihood with respect to matrix $A$ only yields

   $$\hat{A} = \Sigma_\eta^{1/2}Q(\Lambda_K - I_K)^{1/2},$$

where $\Lambda_K$ is the diagonal matrix with the $K$ largest eigenvalue of the matrix $\Sigma_\eta^{-1/2}\tilde{\Gamma}_y(0)\Sigma_\eta^{-1/2}$ on its diagonal and $Q$ is the matrix whose columns are the corresponding eigenvectors. On the other hand, given $A$, the likelihood is maximized by letting $\Sigma_\eta = \text{diag}(\tilde{\Gamma}_y(0) - AA')$. By substituting $\hat{A}$ to $A$ we have the approximate formula

$$\hat{\Sigma}_\eta = \text{diag}(\tilde{\Gamma}_y(0) - \hat{A}\hat{A}').$$

Given an appropriate starting value for $\Sigma_\eta$ we may apply iteratively the formulas that yield $\hat{A}$ and $\hat{\Sigma}_\eta$ until some convergence criterion is met. For instance, the algorithm may stop if the difference between two consecutive values of the maximized log-likelihood is less than a pre-specified tolerance constant. Proofs of formulas are provided by (13), pp. 367-370, who warn that this method does not guarantee convergence. Nevertheless the algorithm is simple and effective in most cases, and the entries on the diagonal of $\Sigma_\eta$ are not constrained to be all equal.

6. Finally the outlier size $\omega$ may be estimated. Note that $\hat{x}_{t_0}$ may possibly include a within-factor outlier $\alpha$. For the static model, if an outlier has been detected at $t = t_0$ and given $\hat{A}$ and $\hat{X}$, we may write the log-likelihood of $\eta_{t_0}$

$$L = -\frac{N}{2}\log(2\pi) - \frac{N}{2}\log|\hat{A}\hat{A}' + \Sigma_\eta| - \frac{1}{2}(y_{t_0} - \hat{A}\hat{x}_{t_0} - \zeta)'(\hat{A}\hat{A}' + \Sigma_\eta)^{-1}(y_{t_0} - \hat{A}\hat{x}_{t_0} - \zeta).$$

Maximization of the likelihood with respect to $\zeta$ yields

$$\hat{\zeta} = y_{t_0} - \hat{A}\hat{x}_{t_0}.$$

To estimate the outlier size $\omega$ we will obtain an estimate $\hat{\alpha}$ so that

$$\hat{\omega} = \hat{A}\hat{\alpha} + \hat{\zeta}. \tag{7.1}$$

According to (7.1) the vector $\hat{\omega}$ may be written as the sum of a vector that is obtained as a linear combination of the columns of $\hat{A}$ and a vector which is orthogonal to the space spanned by the column of $\hat{A}$. The outlier $\alpha$ impacts the factors while $\zeta$ impacts the model structure as a whole. We may estimate the size of $\alpha$ for each one of the components $\hat{x}_1,\ldots,\hat{x}_K$ essentially by building the linear interpolator of each $\hat{x}_{it_0}$, $i = 1,\ldots,K$, based on values observed at $t \neq t_0$, and identifying an outlier at each time when the interpolator is too different from the estimated $\hat{x}_{it_0}$ (see e.g. 4, for details).

As an example, let us simulate $T = 100$ observations of the multivariate $N \times 1$ time series $\{y_t\}$ with $N = 5$ generated by the model

$$y_t = Ax_t + \omega\Delta_t + \eta_t.$$

We assume $K = 4$ factors, so that the matrix $A$ has 5 rows and 4 columns. Let

$$A = \begin{bmatrix} 1 & 0 & 0 & 0 \\ .5 & 1 & 0 & 0 \\ 0 & .5 & 1 & 0 \\ 0 & 0 & .5 & 1 \\ 0 & 0 & 0 & .5 \end{bmatrix}.$$

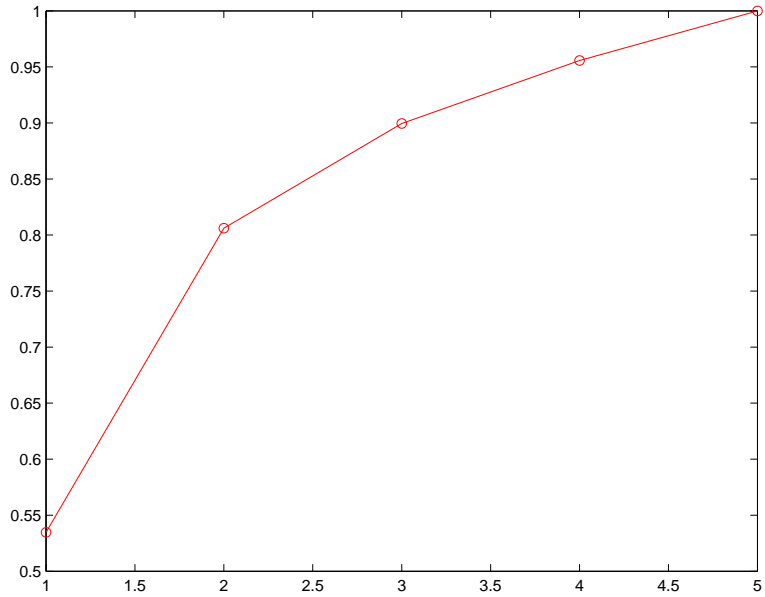

FIG 1. *Plot of eigenvalues of the observed time series spectral matrix at frequency zero*

The factors have been simulated according to the VAR model

$$x_t = \Phi x_{t-1} + \varepsilon_t,$$

where $\Phi = \text{diag}(.7, -.7, .7, -.7)$ and $\{\varepsilon_t\}$ is normal white noise with zero mean and variance $\sigma_\varepsilon^2 = 1 - 0.7^2$. This ensures that the (theoretical) variance of each factor is unity. In addition, as $\Phi$ is diagonal, and by normality assumption, (theoretical) factors are independent. The idiosyncratic component $\{\eta_t\}$ is assumed a zero mean normal white noise sequence with variance $\sigma_\eta^2 = 0.04$. The outlier was located at $t = 100$, that is at the end of the series. Most outlier detection methods are not able to discover potential outliers at the end (or beginning) of the observed time series.

Outlier size was $\omega = (1.5, -1, 0, -4, 5)'$. Each component of the generated time series $y_t$ has (theoretical) variance equal to 1.29, excepted the first and the last one that have variances 1.04 and 0.29 respectively. The outlying observation is rather large only compared to component series 4 and 5, in the remaining cases the outlier size does not exceed twice the standard error of the component series. The standard errors of the simulated time series $\{y_t\}$ with outlier are $(1.0970, 1.2858, 1.1967, 1.1906, 0.7415)$.

The assessment of the number of factors has been performed by examining the eigenvalues of both variance-covariance and spectral density matrices. In Fig. 1 the cumulated eigenvalues of $\hat{F}(0)$ are plotted. It has to be noticed that the smallest eigenvalue is greater than zero, so that the threshold that serves as a decision rule about the number of factors has been computed accordingly (see point 4 above in this Section). This way the correct number of factors $K = 4$ may be identified.

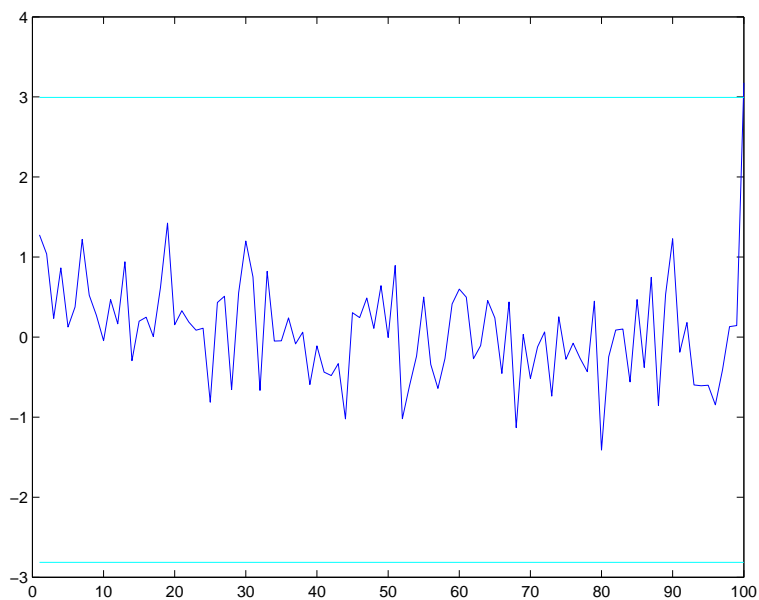

FIG 2. *Univariate time series obtained by projecting the multivariate one along optimal direction*

The optimal direction for detecting outlying observations has been computed

$$z = (-0.5684, 0.4436, -0.3448, -0.1664, 0.5775)'.$$

The univariate time series obtained as the linear combination $\{z'y_t\}$ is displayed in Fig. 2. The outlier at $t = 100$ is clearly highlighted. Mean and standard error of $\{z'y_t\}$ have been computed equal to 0.0902 and 0.6463 respectively. The standardized value at $t_0 = 100$ results equal to 4.7738, larger than the Tchebychev upper bound 4.47 which corresponds to the 5% level.

A VAR has been estimated for the observed time series $\{y_t\}$ and the one-step-ahead forecast has been taken to replace the last observation. The corrected time series was then used to compute the variance-covariance matrix, the covariance matrices at lags 1—4 and the spectral density matrices for 100 Fourier frequencies from $-\pi$ to $\pi$. For checking that the estimated covariance matrices could be assumed symmetric, the interval $(0, \pi)$ has been divided in 4 non overlapping intervals, each of which included 25 frequencies. We obtained for the test statistic the values 10.64, 8.93, 8.77 and 12.32, with 25 degrees of freedom: the critical value at the 5% level is 37.65. As it is greater than the computed statistics, we may not reject the dynamic factor model hypothesis.

The estimates $\hat{A}$ and $\hat{X}$ have been computed by using the three techniques described in this Section, that is approximate temporal decorrelation, eigenvalues-based and iterative maximization of the likelihood function. The latter two methods are similar and yield indeed similar results. The first method is an entirely different approach that takes explicitly into account the covariance matrices at higher lags.

Nevertheless, the estimated dynamic factor model fits the data quite well no matter what method has been used. We report only the plot of the simulated time series (blue line) along with its estimate yielded by the approximate temporal decorrelation algorithm (red line). Other methods yield estimates that overlap almost exactly. In Fig. 3 the observed (simulated) and forecasted (estimated) series are plotted for each component. We may notice that the outlier is not generally apparent by the visual inspection of the graphics.

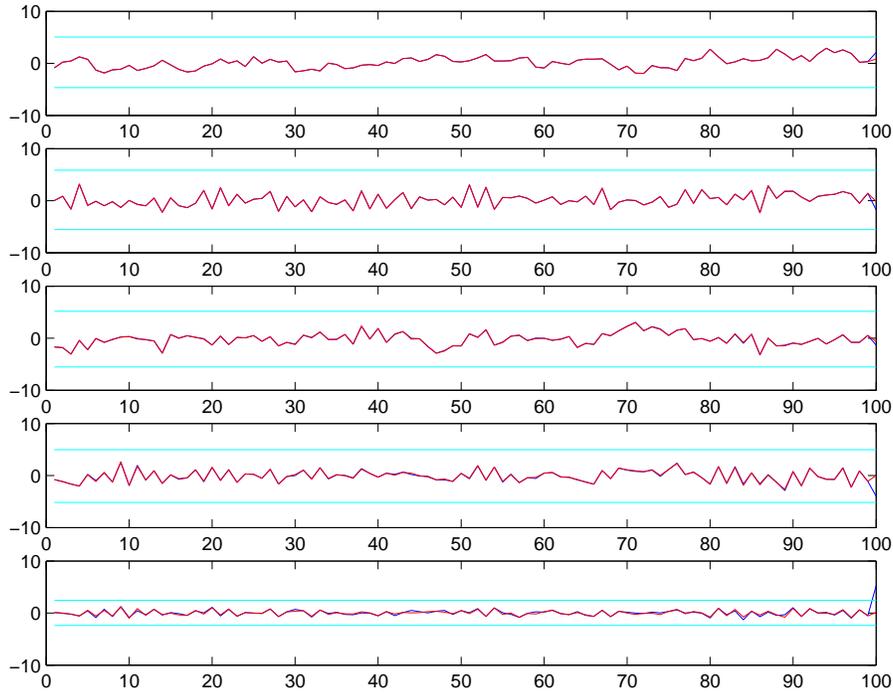

FIG 3. *Observed time series (blue line) and estimated dynamic factor model (red line)*

Then the outlier size has been estimated as the sum of the two components, the first one in the dynamic factor model and the second one in the factors. The two components are orthogonal to each other. Also the first one is orthogonal to the space spanned by the columns of $\hat{A}$ while the other one impacts the dynamic factor model as a linear combination of the columns of $\hat{A}$. The estimates are displayed in Table 1. The three methods yield similar estimates of the total outlier size $\hat{\omega}$ and of the component $\hat{\zeta}$ that impacts the overall model. The sizes of the outlier $\hat{\alpha}$ within the factors differ because the estimated factors themselves depend on the matrices $\hat{A}$ estimated by each of the three methods. The differences are small, however, if we compare the arrays $\hat{A}\hat{\alpha}$.

The results reported in Table 1 seem reliable as regards the recovering of the outlier size $\omega$. We may compute from the 'true' outlier size $\omega$ the arrays

$$\alpha = (1.161, -0.903, -0.903, -0.839)'$$

TABLE 1
*Outlier size estimates yielded by temporal decorrelation, eigenvalue analysis and maximum likelihood methods*

|  | temporal decorrelation | eigenvalue analysis | maximum likelihood |
|---|---|---|---|
| $\hat{\zeta}$ | 0.0944 | 0.2326 | 0.6805 |
|  | −0.0534 | −0.6887 | −1.1582 |
|  | 1.2871 | 1.4950 | 1.9913 |
|  | −2.8438 | −2.7532 | −2.8597 |
|  | 5.4193 | 5.3500 | 4.9343 |
| $\hat{\alpha}$ | 0.4381 | 1.4664 | 0.9291 |
|  | −1.5874 | 1.1927 | 1.8030 |
|  | 1.2749 | −0.4097 | 0.9445 |
|  | −0.1808 | 0.0576 | 1.2409 |
| $\hat{A}\hat{\alpha}$ | 1.1473 | 1.0374 | 0.6822 |
|  | −1.0463 | −0.9231 | −0.5516 |
|  | −2.0511 | −2.2148 | −2.6060 |
|  | −1.1319 | −1.2080 | −1.1289 |
|  | −0.2297 | −0.1667 | −0.1739 |
| $\hat{\omega}$ | 1.2417 | 1.2701 | 1.3627 |
|  | −1.5793 | −1.6118 | −1.7098 |
|  | −0.7640 | −0.7198 | −0.6147 |
|  | −3.9757 | −3.9613 | −3.9886 |
|  | 5.1896 | 5.1833 | 5.1082 |

and

$$\zeta = (0.3387, -0.6774, 1.3548, -2.7097, 5.4194)'.$$

This latter is close to its estimated counterpart (in each of the three versions). As far as $\alpha$ is concerned we have to consider that the product

$$A\alpha = (1.1613, -0.3226, -1.3548, -1.2903, -0.4194)'$$

is close to $\hat{A}\hat{\alpha}$.

## 8. A simulation experiment

We performed a simulation experiment by replicating 1000 times a dynamic factor model and applying three methods for outlier detection and estimation. This way we wanted to test the effectiveness of the method that we are proposing in this paper (let us call it ODFM). Then, we made a comparison with two methods that were available for detecting and estimating outliers in multivariate time series. The first one was proposed by (23) to detect and estimate four types of outliers in multivariate time series modeled by a vector autoregressive integrated moving-average (VARIMA) model (let us call it OARMA). The second one was proposed by (8) as a projection pursuit approach to detect and estimate four types of outliers in multivariate time series not necessarily generated by a VARIMA model (let call it OPP).

We confine our attention only to the most common types of outliers, namely the additive outliers (AO). An AO impacts the series only at the time of its occurrence, while neighboring observations remain unaffected. Other outlier types were defined in the literature: innovation outliers (IO), level changes (LS) and temporary changes (TC). However, IO's are only defined when the data are assumed to follow a VARMA model, which is not our case. LS's arise when the mean levels of each component series change at once, and then remain constant. They are equivalent to AO in the difference, and may be identified by analyzing the differenced data. Finally, a TC in multivariate time series data is defined at $t = t_0$ if a constant $\omega$ which defines the outlier size is added to $y_{t_0}$ and $\delta^k \omega$ is added to $y_{t_0+k}$, $k > 0$, where $0 < \delta < 1$ is a scalar constant. We feel that a TC is a very unlikely behavior in real data, in any case it is easily identified by the existence of an exponentially decaying impulse at rate $\delta$ in the univariate projection series $z'y_t$.

The method of (23) assumes that the ($N$-dimensional) multivariate time series $\{y_t\}$ may be modeled as
$$y_t = x_t + \alpha(B)\omega \xi_t^{(h)},$$
where the unobservable multivariate time series $\{x_t\}$ is generated by the VARIMA model
$$\Phi(B)x_t = c + \Theta(B)\varepsilon_t.$$
In the latter equality,
$$\Phi(B) = I - \Phi_1 B - \ldots - \Phi_p B^p$$
and
$$\Theta(B) = I - \Theta_1 B - \ldots - \Theta_q B^q$$
are $N \times N$ matrix polynomials of finite degrees $p$ and $q$, $c$ is a $N$-dimensional constant vector, and $\{\varepsilon_t\}$ is a sequence of independent and identically distributed normal random vectors with zero mean and covariance matrix $\Sigma_\varepsilon$. Some assumptions are needed to ensure that
$$x_t = c_* + \Phi(B)^{-1}\Theta(B)\varepsilon_t = c_* + \Psi(B)\varepsilon_t$$
is a well defined moving average model. Then, $\alpha(B) = \Psi(B)$ defines an IO and $\alpha(B) = I$ an AO. The date of the outlier is defined by the binary variable $\xi_t^{(h)}$ which equals 1 if $t = h$ (that is, the outlier occurs at $t = h$) and 0 otherwise.

The method of (8) aims at discovering the univariate projections of the multivariate time series that best highlight the presence of outliers. The directions that yield the most useful projections are given by the projections that either maximize or minimize the kurtosis. Moreover, orthogonal directions are to be taken into account as well. The number of projections to be examined for outlier detection is $2N$, where $N$ denotes the dimension of the multivariate time series. If IO's have to be detected, then the method applies to the residual multivariate time series computed from a suitable model fitted to the observed data.

We simulated 200 observations from the $N$-dimensional dynamic factor model
$$y_t = Ax_t + \omega \Delta_t + \eta_t.$$
We assumed $N = 20$ and $K = 4$ factors, simulated according to the VARMA model
$$x_t - \Phi x_{t-1} = \varepsilon_t - \Theta \varepsilon_{t-1}, \tag{8.1}$$

where $\Phi = \text{diag}(0.7, -0.5, 0.5, -0.7)$, $\Theta = \text{diag}(-0.5, 0.7, -0.7, 0.5)$ and $\{\varepsilon_t\}$ is normal white noise with zero mean and variance-covariance matrix $\Sigma_\varepsilon = I_K - \Theta\Phi' - (\Phi - \Theta)\Theta'$. This choice of $\Sigma_\varepsilon$ ensures that the (theoretical) variance-covariance matrix of $\{x_t\}$ is the unit matrix. The matrix $A$ was chosen to have 20 rows and 4 columns. We let

$$A = \begin{bmatrix} 2 & 1 & 0 & 0 \\ 1 & 0 & 2 & 0 \\ 1 & 0 & 0 & 2 \\ 0 & 2 & 1 & 0 \\ 0 & 1 & 0 & 2 \\ 0 & 0 & 1 & -2 \\ 1 & 1 & -2 & 0 \\ 1 & -2 & 0 & 1 \\ -2 & 0 & 1 & 1 \\ 0 & 1 & -2 & 1 \\ 1 & 2 & 0 & 0 \\ 2 & 0 & 1 & 0 \\ 2 & 0 & 0 & 1 \\ 0 & 1 & 2 & 0 \\ 0 & 2 & 0 & 1 \\ 0 & 0 & -2 & 1 \\ 1 & -2 & 1 & 0 \\ -2 & 1 & 0 & 1 \\ 1 & 0 & 1 & -2 \\ 0 & -2 & 1 & 1 \end{bmatrix}.$$

The matrix $A$ may be verified to have full rank. The (theoretical) factors are independent because both $\Phi$ and $\Theta$ are diagonal matrices and the $\{\varepsilon_t\}$'s are uncorrelated normal random variables. All idiosyncratic components $\{\eta_t\}$ were assumed zero mean normal with variance $\sigma_\eta^2 = 0.04$.

We checked two outlier configurations. The first one was an isolated multivariate outlier at $t = 100$. The second one was a patch, that is a sequence of neighboring outliers at $t = 99, 100, 101$. The size of each and every outlier was chosen equal to 0.6. This figure was chosen in comparison with the standard error $\sigma_\eta = 0.2$ of each of the idiosyncratic components. The total multivariate outlier size is the $N$-dimensional vector $\omega$ with entries all equal to 0.6. The outlier $\omega$ may be split in a term $\zeta$ orthogonal to the columns of $A$, that is

$$\zeta = (I - A(A'A)^{-1}A')\omega,$$

and a term $A\alpha$ which is a linear combination of the columns of $A$ with coefficients

$$\alpha = (A'A)^{-1}A'\omega.$$

The coefficients $\alpha$ may be thought of as the sizes of an outlier that impacts the factors.

All computations were performed by using the Matlab package. For each of 1000 replications we generated 300 independent identically normally distributed $K$-variate random vectors and 200 independent identically distributed $N$-dimensional random vectors all with mean zero and variance-covariance equal to the unit matrix. From the 300 $K$-dimensional random vectors 300 observations were generated from the ARMA model (8.1). Then, the data were transformed as explained before, to obtain unit variance factors. The first 100 $K$-dimensional data were discarded to remove the effect of the (random) initial values. As a result, a $K$-dimensional factor time series of length 200 was obtained. The $N$-dimensional white noise was pre-multiplied by the inverse of the square root of the matrix $\Sigma_\eta$. This was the artificial idiosyncratic component that was added to the factor data. Finally the two outlier structures were superimposed to the artificial data generated from the dynamic factor model. In each replication, and for each of the two outlier structures, the methods ODFM, OARMA and OPP were applied for outlier detection and estimation. The usual Monte Carlo simulation procedures were used to compute the percentages of both correct and false identifications and the average and standard errors of the estimates. A synthetic measure of the distance between the estimated and true outlier was obtained by computing the norm of the vector difference between the estimated and true outlier size.

TABLE 2
*Percentages of estimated number of factors in the presence of outliers*

| Outlier type | Estimated number of factors | |
|---|---|---|
| | $K=3$ | $K=4$ |
| isolated $t=100$ | 0.3% | 99.7% |
| patch $t=99,100,101$ | 0.2% | 99.8% |

In the present context it seems of interest to report some results concerned with the validity of a dynamic factor model to fit the data. We divided the frequency interval $(0,\pi)$ into four sub-intervals of equal size, namely $[\frac{2\pi}{T},\frac{\pi}{4}]$, $[\frac{\pi}{4}+\frac{2\pi}{T},\frac{\pi}{2}]$, $[\frac{\pi}{2}+\frac{2\pi}{T},\frac{3}{4}\pi]$, and $[\frac{3}{4}\pi+\frac{2\pi}{T},\pi]$. Each sub-interval included $T/8-1$ frequencies (here 24 frequencies as $T=200$), and the sub-interval centers were $\frac{\pi}{8}$, $\frac{3}{8}\pi$, $\frac{5}{8}\pi$, and $\frac{7}{8}\pi$ respectively. The null hypothesis, $H_0$: variance-covariance matrices are symmetric, tested using the LRT statistics of Theorem 4.2 was never rejected at significance level 5% neither in the presence of an isolated outlier nor an outlier patch.

Table 2 shows that the number of factors $(K=4)$ was correctly estimated in almost all replications.

For method ODFM we computed the $N-K$ univariate projections obtained as linear combination of the multivariate time series

$$w_{i,t} = (v_i)'y_t, \quad i=1,\ldots,N-K,$$

where $v_i$ is the eigenvector of the variance-covariance matrix $\hat{\Gamma}_y(0)$ of the observed multivariate time series associated with the eigenvalue $\lambda_i$. The eigenvalues were arranged in ascending order, that is the eigenvectors $v_1,\ldots,v_{N-K}$ belong to the smallest eigenvalues $\lambda_1,\ldots,\lambda_{N-K}$ respectively. Then, the presence of an outlier in the multivariate time series $\{y_t\}$ was detected at time $t$ if

$$|w_{i,t} - \bar{w}_i| > k_\alpha \sigma_{w_i}$$

for some (that is, at least one) 1-dimensional time series $\{w_{i,t}\}$. $\bar{w}_i$ and $\sigma_{w_i}$ were the $w_{i,t}$ sample mean and standard deviation respectively. The threshold parameter $k_\alpha$ may be computed according to the Tchebychev inequality. We chose the significance level $\alpha = 0.05$ so that approximately $k_\alpha = 4.47$.

The parameters of the dynamic factor model were estimated along the guidelines given by Theorem 5.1 in Section 5 (see step 5(b) in Section 7 as well). The difference between the observed time series and the estimated dynamic factor model values was assumed to yield the estimate of the outlier size $\hat{\zeta}$. Then outliers $\hat{\alpha}$ were estimated in each factor by using Formula (2.2b) p. 194 of (4). The total outlier size was obtained as $\hat{\omega} = \hat{\zeta} + \hat{A}\hat{\alpha}$.

The OARMA method was implemented along the guidelines given by (23). A VAR model of order $M = 4$ was fitted to the multivariate time series $\{y_t\}$. We assumed that only outliers of either additive or innovation type could be present. The Mahalanobis type statistic for either type of outliers

$$J_{i,h} = (\hat{\omega}_{i,h})' \Sigma_{i,h}^{-1} \hat{\omega}_{i,h}$$

was computed for each time $t = h$ and $i = I$ for IO and $i = A$ for AO. $\Sigma_{i,h}$ denotes the covariance matrix of the estimator. If the maximum across time of $J_{A,h}$ or $J_{I,h}$ exceeded their respective 95-th percentile, then either an AO or an IO was assumed at $h = h_{\max}$ according to which statistics $J_{\max}(i, h_i) = \max_h J_{i,h}$, $i = I, A$, was the greatest. Then the outlier size was estimated and its effect removed from the multivariate time series. The procedure was iterated until no more outliers were found. Tables of percentiles of the statistics $J_{\max}(I, h_I)$ or $J_{\max}(A, h_A)$ are available only up to dimension 10 (see (23), Table 1 p. 797, and (8), Table 4 p. 664). So we computed empirical percentiles from 10000 artificial multivariate time series generated by the dynamic factor model with $N = 20$ and $T = 200$. We obtained $J_{\max}(A, h_A) = 47.7911$ and $J_{\max}(I, h_I) = 46.6493$ at the 5% significance level.

The estimates of the outlier size were computed by using the Formulas $\hat{\omega}_{I,h}$ for innovation outliers and $\hat{\omega}_{A,h}$ for additive outliers provided by (23) p. 794.

As far as the OPP method is concerned the direction that maximizes the kurtosis of the linear projection of the multivariate time series $\{y_t\}$ and all orthogonal directions had to be computed. The direction that minimizes the kurtosis and its orthogonal directions had to be computed as well. To compute these $2N$ projections we used the Matlab routines by (15) available on the web (http:// halweb.uc3m.es/fjp/download.html). In this case too we confined ourselves only to AO and IO. In this latter case, the procedure was applied to the multivariate residual time series $\{a_t\}$ computed by fitting a VAR of order $M = 4$ to $\{y_t\}$. Then, each and every projection was searched for outliers by using the univariate counterpart of the statistic in the OARMA method. The maximum across time and across projections was computed and let $\Lambda_A$ denote the maximum found on the projections of the multivariate time series $\{y_t\}$ and $\Lambda_I$ denote the maximum found on the projections of the multivariate residual time series $\{a_t\}$. Both $\Lambda_A$ and $\Lambda_I$ were compared with their appropriate thresholds and either an AO or an IO was detected if the greatest of $\Lambda_A$ and $\Lambda_I$ exceeded its threshold. As tables of percentiles of $\Lambda_A$ and $\Lambda_I$ are available only up to $N = 10$ (see (8), Table 2, p. 663), we computed the 95-th percentiles by simulation from the dynamic factor model with $N = 20$ and $T = 200$. In this case we obtained after 10000 replications $\Lambda_A = 5.9772$ and $\Lambda_I = 6.0392$.

TABLE 3
*Isolated outlier percentage detection by the ODFM, OARMA and OPP methods*

| ODFM outlier id. percent | | OARMA outlier id. percentage | | | OPP outlier id. percentage | | |
|---|---|---|---|---|---|---|---|
| correct | false | correct AO | correct IO | false | correct AO | correct IO | false |
| 96.9 | 1.7 | 2.56 | 94.64 | 100 | 58.2 | 40.9 | 4.3 |
| (3.5199) | (2.3473) | (2.21) | (2.99) | | (8.1462) | (8.2698) | (4.2320) |

Once date and type of outliers were determined, to estimate their size we fitted a VAR model of order $M = 4$ to the multivariate time series $\{y_t\}$ and computed $\hat{\omega}_{I,h}$ for innovation outliers and $\hat{\omega}_{A,h}$ for additive outliers as in (23) p. 794.

The results obtained by applying the three procedures to discover the outlying observation in the artificial data set with an isolated outlier are displayed in Table 3. The Monte Carlo statistics were computed on 1000 replications. In order to compute the standard errors of the estimates of the percentages of detection, we divided the replications in 40 groups of length 25 each. The percentages were computed for each group of replications, then we could compute the standard error of each percentage by using 40 values. The percentage of correct detection of the isolated outlier in $t = 100$ was greater than 95% for all methods. The standard errors are comparable as regards methods ODFM and OARMA but the standard error is slightly larger for the OPP method. Both OARMA and OPP methods identified in most cases an outlier of innovation type possibly because the IO allows for greater flexibility as far as fitting the observed data is concerned. On the other hand, if we constrained the OARMA and OPP methods to search for AO only, then the percentages of correct detection dropped dramatically. We recorded false outlier detections as well. The number of replications (percentage) where one or more wrong detections occurred was low for methods ODFM and OPP while the method OARMA wrongly detected outlying observations in every replication. As overall figures, the observations detected as outlying ones over 1000 replications (series length 200) were 17, 80, and 9003 for the ODFM, OPP, and OARMA methods respectively.

The results for the data set with a patch at times $t = 99, 100, 101$ are displayed in Table 4. Considering the outlier separately, the ODFM method detected each of the three outliers about 99%, while for the OARMA and OPP methods high percentages were recorded only for the detection of the first outlier in the patch. Again, these latter methods almost always identified the outlying observations as IO. As a consequence, the outliers in $t = 100$ and $t = 101$ could be well explained by the innovation outlier structure. The percentages of false identifications were rather low for the ODFM and OPP methods, while the OARMA method wrongly identified outliers in all 1000 multivariate time series. The average number of false outliers was in this case about 7, that is less than in the case of multivariate time series with an isolated outlier. By considering the patch as a whole, correct identification was performed 97.5% by ODFM, only about 58% by OARMA and never by OPP. At least an outlier in the patch was detected by ODFM, OARMA and OPP 100%, 85.09% and 85.4% respectively.

Outlier size estimates are displayed in Table 5 for the case of the isolated outlier and in Table 6 for the case of the outlier patch. The adequacy of the estimates was evaluated by the norm of the difference between the estimated and true outlier size vectors. This figure yields a kind of distance measure which accounts both for bias and variability of the estimates. We may remind that the outlier size is equal to 0.6 for all component time series. For the ODFM method we distinguish the total outlier from its orthogonal part. This is of interest because the dynamic factor model structure is directly impacted by the array orthogonal to the columns of the matrix *A*, while the rest impacts the factors. For the other methods, OARMA and OPP, we distinguish between AO and IO identification. We may notice that as far as the ODFM method is concerned, the orthogonal part is close to its true counterpart while the estimated total outlier seems more variable. The distance is about 7 in all cases. As regards the other methods OARMA and OPP the distance between the estimated and true outlier size is smaller, approximately between 4 and 5. In the case of OPP this figure is even smaller (about 2.5) if the outlier is identified as AO.

TABLE 4
*Outlier patch percentage detection by the ODFM, OARMA and OPP methods*

| Method | ODFM id. percent | OARMA iden. percentage | | OPP iden. percentage | |
|---|---|---|---|---|---|
| | | AO | IO | AO | IO |
| Outlier $t = 99$ | 98.9 (1.9975) | 0.1 (0.24) | 83.09 (3.79) | 0.0 | 82.9 (7.9492) |
| Outlier $t = 100$ | 99.8 (0.8718) | 0.2 (0.42) | 63.72 (4.79) | 0.0 | 1.4 (2.1071) |
| Outlier $t = 101$ | 98.8 (2.7129) | 4.34 (2.23) | 60.06 (4.69) | 2.9 (3.3451) | 0.0 |
| False detection | 1.5 (2.1331) | 100 | | 4.8 (3.9192) | |

Such distances depend both on bias and variability. To distinguish between these two sources, we present the bias and standard errors of the estimates of $\omega_i$ in Table 7 for the case of isolated outlier. Figures are bias and standard errors of estimated outlier sizes in each component, over the correctly identified replications, and averaged on the 20 components. It may be seen that the proposed ODFM method provides less biased estimates of the outlier size, while the variability is larger than the projection pursuit method. Since the distances of the estimated $\hat{\zeta}$ from the true values are generally small (see Table 5 and 6), we conclude that estimation of univariate outliers in the (estimated) factors, or interaction between them and the estimates of the factor matrix *A*, are responsible for such a larger variability. This may be caused by the method employed in Theorem 5.1 for estimating $\{x_t\}$. A more efficient method is currently under investigation.

## 9. Application to real data

We used a quarterly economic data set which covers three countries, France, Germany and Italy, from 1960 to 1999, for illustrating method ODFM. The data are taken from the seven-country data set used by (22) to compute combination forecasts of output growth. The preliminary transformations suggested in the aforementioned paper were applied. Then, 154 transformed data for each time series were available from III quarter 1960 to IV quarter 1999. We considered only the time series indicated by 'a' in Table Ib in (22).

TABLE 5
*Deviation from true values of isolated outlier size estimates yielded by the ODFM, OARMA and OPP methods*

| Method | ODFM | | OARMA | | OPP | |
|---|---|---|---|---|---|---|
| | $\|\hat{\zeta} - \zeta\|$ | $\|\hat{\omega} - \omega\|$ | $\|\hat{\omega}_A - \omega\|$ | $\|\hat{\omega}_I - \omega\|$ | $\|\hat{\omega}_A - \omega\|$ | $\|\hat{\omega}_I - \omega\|$ |
| Outlier | 0.7805 | 7.06 | 3.3778 | 4.5135 | 2.5843 | 4.3167 |
| $t = 100$ | (0.1364) | (2.8675) | (1.4531) | (3.4293) | (0.6581) | (1.5351) |

TABLE 6
*Deviation from true values of outlier patch size estimates yielded by the ODFM, OARMA and OPP methods*

| Method | ODFM | | OARMA | | OPP | |
|---|---|---|---|---|---|---|
| | $\|\hat{\zeta} - \zeta\|$ | $\|\hat{\omega} - \omega\|$ | $\|\hat{\omega}_A - \omega\|$ | $\|\hat{\omega}_I - \omega\|$ | $\|\hat{\omega}_A - \omega\|$ | $\|\hat{\omega}_I - \omega\|$ |
| Outlier | 0.7817 | 7.1861 | 3.9729 | 4.5533 | | 4.3767 |
| $t = 99$ | (0.1386) | (2.8408) | (1.29) | (1.6458) | | (1.5117) |
| Outlier | 0.7829 | 7.1386 | 3.9406 | 6.1979 | | 4.6051 |
| $t = 100$ | (0.1381) | (2.8656) | (1.1391) | (2.4194) | | (1.0134) |
| Outlier | 0.7773 | 7.2700 | 4.0333 | 5.1087 | 2.4556 | |
| $t = 101$ | (0.1415) | (2.7589) | (1.4710) | (1.9191) | (0.4919) | |

TABLE 7
*Bias and standard errors of the outlier size estimates, averaged on the 20 components (isolated outlier at $t = 100$)*

| Method | ODFM | OPP–AO | OPP–IO |
|---|---|---|---|
| correctly identified repl.ns | 969 | 582 | 409 |
| average bias | 0.014 | -0.373 | -0.231 |
| average standard error | 1.70 | 0.46 | 0.99 |

TABLE 8
*French (f), German (g) and Italian (i) quarterly economic data and preliminary transformations*

| Label | Description | Transform. | Countries | | |
|---|---|---|---|---|---|
| cpi | Consumer price index | $\Delta^2 \ln$ | f | g | i |
| ip | Index of industrial production | $\Delta \ln$ | f | g | i |
| pgdp | Gross domestic product deflator | $\Delta^2 \ln$ | | g | i |
| rbndl | Interest rate of long-term government bonds | $\Delta$ | f | g | i |
| rbndm | Interest rate of medium-term government bonds | $\Delta$ | | | i |
| rcommod | Real commodity price index | $\Delta \ln$ | f | g | i |
| rgdp | Real gross domestic product | $\Delta \ln$ | | g | i |
| roil | Real oil prices | $\Delta \ln$ | f | g | i |
| rgold | Real gold prices | $\Delta \ln$ | f | g | i |
| rstock | Real stock price index | $\Delta \ln$ | f | g | i |
| unemp | Unemployment rate | $\Delta$ | | | i |

For France 7 series were available from 1960 to 1999, 9 series for Germany and 11 series for Italy. Time series labels, description, transformations and countries are displayed in Table 8.

The data set was composed of 27 time series of length $T = 154$. The number of factors was checked by using the eigenvalues of the variance-covariance matrix. The eigenvalues were computed and arranged in descending order. The cumulated sum is plotted in Fig. 4. The smallest integer such that the cumulated sum exceeded 0.95 was 6, so that we assumed the number of factor $K = 6$.

The symmetry of the variance-covariance matrices computed for lags 1, 2, 3 and 4 was checked along the guidelines given in Section 4. Two statistics for this test were computed by averaging the periodogram in two frequency bands, centered in $\lambda = \pi/4$ the first one and in $\lambda = 3\pi/4$ the second one. We obtained respectively 174.81 and 143.88 which are not significant, so that the null hypothesis was not rejected.

The identification stage required the computation of $N - K = 14$ projections that were searched for outlying observations along the guidelines given in Section 8. The projections that exceeded the thresholds computed by the Tchebychev inequality suggested the presence of outliers at times $t = 23, 24, 37, 63, 124$. The graphical display of these projections in Fig. 5 shows that each projection disclosed an outlier separately, excepted the third projection (outliers at $t = 23$ and $t = 24$).

Examination of the estimated outlier sizes allowed the time series that most contributed to the multivariate outlier to be discovered. We could see that the outlier at $t = 23$ (I/1966) was apparent in time series rbndl (g) and pgdp, rbndl, rstock, unemp (i). The outlier at $t = 24$ (II/1966) was mainly determined by roil (f), roil and rstock (g), and pgdp, rbndl, rstock, unemp (i). The outlier at $t = 37$ (III/1969) was mainly determined by rbndl (f), rbndl (g), and rbndl and unemp (i). The main influence on the outlier at $t = 63$ (I/76) came from rbndl and rgold (f), rbndl and rgold (g), and rbndl, rbndm, rcommod, roil (i). The outlier at $t = 124$ (II/1991) was apparent only in time series rbndl (f). The interest rate of long-term government bonds (rbndl) was present in all outliers, that is this time series produced most of the unexpected values observed in the multivariate time series. Its effect was common to the three countries at $t = 37$ and $t = 63$ while its effect was limited to Germany and Italy at $t = 23$, Italy at $t = 24$, and France at $t = 124$. Some time series were available for Italy only, so it seemed of interest to apply the procedure to the data set which included only the quarterly economic Italian data. We could fit a dynamic factor model with 3 factors and 4 outliers were found, at $t = 23$ (pgdp, rbndl, rbndm, rstock), $t = 24$ (pgdp, rbndl), $t = 55$ (rbndm, roil, rgold, rstock, unemp), and $t = 103$ (rbndl, rbndm, roil, rstock, unemp). Observation at $t = 55$ corresponds to I/1974 and $t = 103$ to I/1986. By comparing the results for Italy with those for the three countries together we could argue that some outliers were 'international' while others regarded only one or two countries. So we had to consider outliers at $t = 23$ and $t = 24$ as 'international' while outliers at $t = 55$ and $t = 103$ as 'national'. Some of the time series that originated these two outliers, that is rbndm and unemp, were present in the Italian data set only. Similarly we could not consider the outlier at $t = 124$ as 'international' but 'national' limited to the quarterly economic French data. Note that the outlier at $t = 23$ was not considered present in the French data, but dates $t = 23$ and $t = 24$ are close so that the outlying observations were possibly related to some common circumstance.

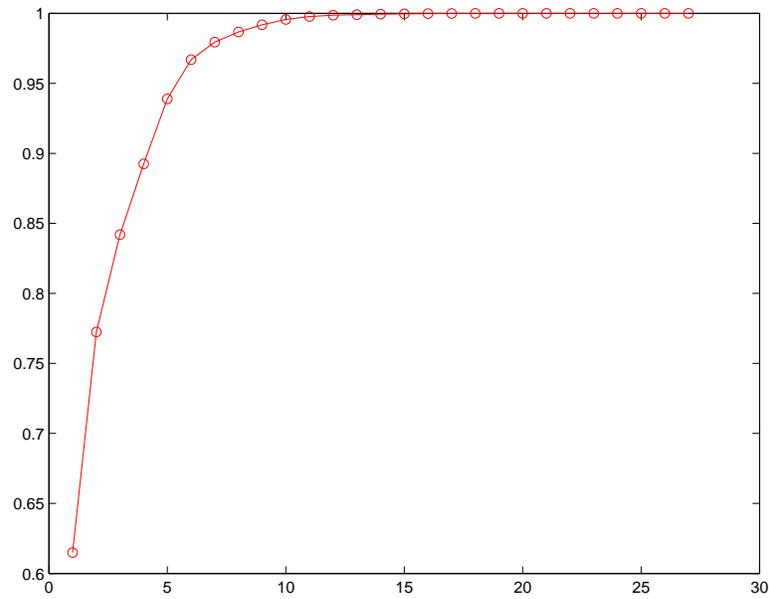

FIG 4. *Cumulated eigenvalues of the variance-covariance matrix of* 27 *economic quarterly data time series*

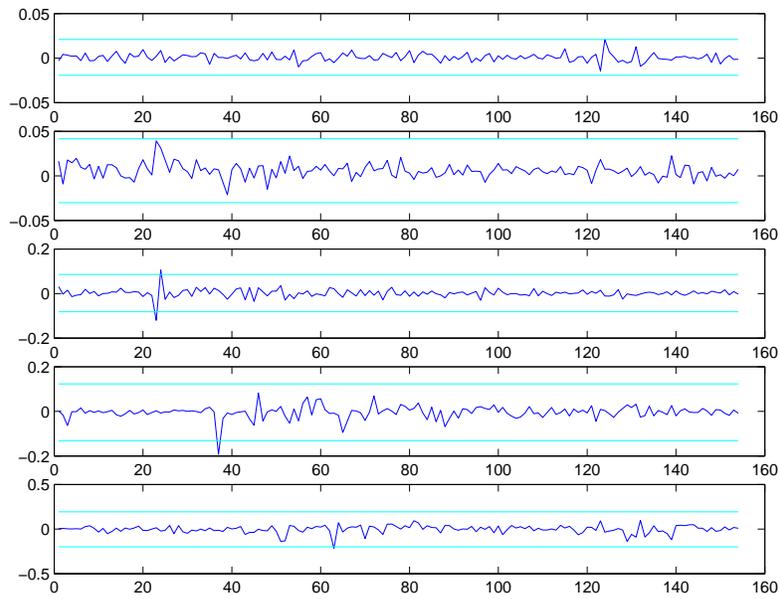

FIG 5. *Projections of the multivariate time series where the presence of outliers may be detected*

## 10. Conclusions

We presented a method to discover outliers in multivariate time series generated by a dynamic factor model. This method was found to yield best results compared to two other methods aimed at discovering outliers in multivariate time series in a different framework. Our method relies on the assumption that the multivariate time series is generated by a dynamic factor model, therefore the procedure to check the dynamic model adequacy for fitting the data should be carefully applied to ensure that genuine outliers could be discovered. If assumptions were carefully checked and requirements met, both the simulation experiment and the application to real data showed that the method presented in this article was effective for outlier detection and estimation, cautious against false outlier identification and, in addition, simple to implement. The estimates of the total outlier size were found less biased, but more variable than those obtained by using the other two methods. On the contrary the estimates of the size of the part of outlier that impacts the dynamic factor model without affecting the factors were found accurate and close to their respective true values. Improvements of the estimation method is currently subject of further research.

**Appendix A: Appendix**

*A.1. Proof of Theorem 2.1*

*Proof.* Suppose that $x_t^* = Cx_t$ and Assumptions 1 and 2 hold, then

$$x_t^* = \sum_{j=0}^{\infty} \Theta_j^* \varepsilon_{t-j},$$

where $\Theta_j^* = C\Theta_j$ and $\Theta_j = \text{diag}(\theta_{11}^{(j)}, \ldots, \theta_{KK}^{(j)})$. But since the $x_t^*$'s are mutually independent, each $\Theta_j^*$, $j = 0, 1, \ldots$, is diagonal, therefore $\sum_k c_{sk} \theta_{kr}^{(j)} = 0$, $r \neq s$. The $\Theta_j$'s are diagonal, therefore $c_{sr}\theta_{rr}^{(j)} = 0$, $r \neq s$. Since for each $r$ at least one $\theta_{rr}^j$ is non-zero (otherwise $x_{r,t}$ would be equal to zero for each and every $t$) it follows that $C$ too has to be diagonal. Thus, again from Assumption 1,

$$\Gamma_x^*(0) = \text{cov}(x_t^*, (x_t^*)') = C\Gamma_x(0)C' = CC' = I,$$

that is $c_{kk}^2 = 1$, $k = 1, \ldots, K$, which implies $C = I$ up to a change of sign. □

*A.2. Proof of Theorem 3.1*

*Proof.* Note first that, for $h \neq 0$, $\Gamma_y(h) = A\Gamma_x(h)A'$ has rank $\leq K$, therefore $\Gamma_y(h)$ has at least $(N-K)$ eigenvalues equal to zero for any $h \neq 0$ (zero is not necessarily the smallest eigenvalue because $\Gamma_y(h)$ is not positive definite).
If $z'A = 0$, then $\Gamma_y(h)z = A\Gamma_x(h)A'z = 0$, therefore $z$ is eigenvector associated with a zero eigenvalue of each matrix $\Gamma_y(h)$ for any $h \neq 0$.
On the other hand, let $z$ be an eigenvector associated with a zero eigenvalue of each matrix $\Gamma_y(h), h \neq 0$ (the set of such vectors is not empty because it includes $V_A^\perp$). Then $\Gamma_y(h)z = A\Gamma_x(h)A'z = 0$ for any $h$ and on multiplying by the generalized inverse $(A'A)^{-1}A'$ we get $\Gamma_x(h)A'z = 0$. If we let $c = A'z$ the last equation reads $\Gamma_x(h)c = 0$ or, since $\Gamma_x(h)$ is diagonal,

$$\gamma_j^x(h)c_j = 0, \quad j = 1, 2, \ldots, K, h = 1, 2, \ldots$$

which under our hypotheses implies $c_j = 0, j = 1, 2, \ldots, K$ and therefore $A'z = 0$. □

*A.3. Proof of Theorem 3.2*

*Proof.* We note first that

$$\Gamma_y(h) = \sum_{u=1}^{s-h} \psi_u \psi_{u+h}', \quad h = 1, 2, \ldots, s-1$$

and $\Gamma_y(h) = 0$ for $h \geq s$. Moreover, the matrix $[\psi_1, \psi_2, \ldots, \psi_s]$ has rank $\leq sK(< N)$, thus there exists a vector $z$ such that $z'\psi_u = 0$ for any $u$. For such a vector we have

$$\Gamma_y(h)z = \sum_{u=1}^{s-h} \psi_u \left( \psi_{u+h}' z \right) = 0, \quad h = 1, 2, \ldots, s-1$$

and it also follows that $\Gamma_y(h)$ has a zero eigenvalue for each $h = 1, 2, \ldots, s-1$.
Suppose now that $z$ is eigenvector associated with a zero eigenvalue of each $\Gamma_y(h)$, and write $c_u = \psi_u' z$. Then
$$0 = \Gamma_y(s-1)z = \psi_1 \psi_s' z = \psi_1 c_s$$
and since the rank of $\psi_1$ is $K$ it follows $c_s = 0$. Now,
$$0 = \Gamma_y(s-2)z = \psi_1 \psi_{s-1}' z + \psi_2 \psi_s' z = \psi_1 c_{s-1} + \psi_2 c_s = \psi_1 c_{s-1}$$
and $c_{s-1} = 0$. On repeating the same argument for lags $s-3, s-2, \ldots, 2, 1$ we obtain $c_{s-2} = c_{s-3} = \ldots = c_2 = 0$. It remains to show that $c_1 = 0$ also. From assumption (ii) for a fixed $k$ the rank of $\psi_k$ is $K$, and
$$0 = z'\Gamma_y(k-1) = z' \sum_{u=1}^{s-k+1} \psi_u \psi_{u+k-1}' = c_1' \psi_k' + c_2' \psi_{k+1}' + c_3' \psi_{k+2}' + \ldots + c_{s-k+1}' \psi_s' = c_1' \psi_k'$$
and it follows $c_1 = 0$, which completes the proof. □

### A.4. Proof of Theorem 3.3

*Proof.* Let us recall that $\operatorname{rank}(AA') = \operatorname{rank}(A) = K$. Then $AA'$ has $N - K$ eigenvalues equal to zero. We have that $\Gamma_y(0) = AA' + \sigma^2 I$ is symmetric and positive definite. From the relationship
$$AA' - \lambda I = \Gamma_y(0) - (\lambda + \sigma^2)I$$
it follows that $\Gamma_y(0)$ has minimum eigenvalues equal to $\sigma^2$ with multiplicity $N - K$. Let $z$ be a corresponding eigenvector. Then,
$$\sigma^2 = z'\Gamma_y(0)z = z'AA'z + \sigma^2,$$
and $\|z'A\| = 0$ follows. On the other hand, if $z'A = 0$ then $\Gamma_y(0)z = AA'z + \sigma^2 z$ and $z$ is eigenvector associated to the smallest eigenvalue $\sigma^2$. □

### A.5. Proof of Theorem 3.4

*Proof.* Since
$$\Gamma_y(0) = \sum_{u=1}^{s} \psi_u \psi_u' + \sigma^2 I$$
and for any $v$ such that $\|v\|^2 = 1$ we have
$$v'\Gamma_y(0)v = \sum_{u=1}^{s} \|v'\psi_u\|^2 + \sigma^2 \geq \sigma^2$$
it follows that the smallest eigenvalue of $\Gamma_y(0)$ is at least $\sigma^2$. If $\psi_u' z = 0, u = 1, 2, \ldots, s$ then
$$\Gamma_y(0)z = \sum_u \psi_u(\psi_u' z) + \sigma^2 z = \sigma^2 z$$

and $z$ is eigenvector associated with $\sigma^2$, the smallest eigenvalue of $\Gamma_y(0)$.

On the other hand, let $z$ denote an eigenvector associated with the smallest eigenvalue $\sigma^2$, it follows that $z'\Gamma_y(0)z = \sigma^2$. Then

$$z'\Gamma_y(0)z = \sum_{u=1}^{s} \|\psi'_u z\|^2 + \sigma^2 = \sigma^2.$$

Hence $\sum_u \|\psi'_u z\|^2 = 0$, therefore $\psi'_u z = 0$ for any $u$. $\square$

### A.6. Proof of Theorem 4.1

*Proof.* If $\Gamma_y(h)$ is symmetric for any $h$ we can write

$$F(\lambda) = \frac{1}{2\pi} \sum_{h=-\infty}^{\infty} \Gamma_y(h)e^{-i\lambda h} = \frac{1}{2\pi}\left\{\Gamma_y(0) + 2\sum_{h=1}^{\infty} \Gamma_y(h)\cos(\lambda h)\right\}$$

which is real for any $\lambda$. On the other side, if $F(\lambda)$ is real then $F(\lambda) = \overline{F(\lambda)} = F(-\lambda)$ and

$$\Gamma_y(h)' = \Gamma_y(-h) = \int F(\lambda)e^{-i\lambda h}d\lambda = \int F(\omega)e^{i\omega h}d\omega = \Gamma_y(h).$$

$\square$

### A.7. Proof of Theorem 4.2

*Proof.* Assuming normality, $(X_j^R, X_j^I)$ and $(X_k^R, X_k^I)$ for $1 \leq j < k \leq J$ are independent, therefore the unrestricted likelihood may be written

$$L\left(X_1^R, X_1^I, ..., X_J^R, X_J^I\right) =$$
$$\pi^{-NJ} \left|\begin{array}{cc} \operatorname{Re}F & -\operatorname{Im}F \\ \operatorname{Im}F & \operatorname{Re}F \end{array}\right|^{-J/2} \exp\left\{-\sum_{j=1}^{J}((X_j^R)', (X_j^I)')\left(\begin{array}{cc} \operatorname{Re}F & -\operatorname{Im}F \\ \operatorname{Im}F & \operatorname{Re}F \end{array}\right)^{-1}\left(\begin{array}{c} X_j^R \\ X_j^I \end{array}\right)\right\}.$$
(A.1)

Therefore the maximum likelihood estimate of $\left(\begin{array}{cc} \operatorname{Re}F & -\operatorname{Im}F \\ \operatorname{Im}F & \operatorname{Re}F \end{array}\right)$ equals $\left(\begin{array}{cc} S_R & -S_I \\ S_I & S_R \end{array}\right)$ (10) and the unrestricted maximum likelihood is

$$\hat{L} = \pi^{-NJ}\left|\begin{array}{cc} S_R & -S_I \\ S_I & S_R \end{array}\right|^{-J/2} \exp\{-2NJ\}.$$

Under $H_0$, $F = \operatorname{Re}F$, both $X_j^R$ and $X_j^I$ are independently $N(0, \frac{1}{2}F)$, thus

$$L(X_1^R, X_1^I, ..., X_J^R, X_J^I) = \pi^{-NJ}|\operatorname{Re}F|^{-J}\exp\left\{-J\operatorname{tr}\left\{(\operatorname{Re}F)^{-1}S_R\right\}\right\}.$$

On maximizing with respect to $\operatorname{Re} F$, we find that the maximum likelihood estimate of $\operatorname{Re} F$ is $S_R$ and the maximum equals

$$\hat{L}_0 = \pi^{-NJ} |S_R|^{-J} \exp\{-2NJ\} .$$

The likelihood ratio is

$$\frac{\hat{L}_0}{\hat{L}} = |S_R|^{-J} \begin{vmatrix} S_R & -S_I \\ S_I & S_R \end{vmatrix}^{J/2} .$$

Using

$$\begin{vmatrix} S_R & -S_I \\ S_I & S_R \end{vmatrix} = |S_R| |S_R + S_I S_R^{-1} S_I|$$

it follows that the likelihood ratio is a monotonic function of the statistic

$$U = |S_R|^{-1} |S_R + S_I S_R^{-1} S_I| = |I + (S_R^{-1} S_I)^2|$$

and the rejection region is $U < c$. The distribution of $U$ is analyzed in Chapter 8 of (1), and corresponds in his notation to $U_{N,N,J-N-1}$. □

### A.8. Proof of Theorem 5.1

*Proof.* The first part of the theorem follows directly from proposition (4) in (12), p. 46. So we only have to prove that $Y\hat{X}' = W_K \Lambda_K^{1/2}$. Let the singular value decomposition of $Y$ be written as in proposition (1) in (12), p. 60, that is

$$Y = WDV'.$$

The $N \times T$ matrix $D$ may be written

$$D = \begin{bmatrix} \Lambda_K^{1/2} & 0_{K \times (r-K)} & 0_{K \times (T-r)} \\ 0_{(r-K) \times K} & \Lambda_{r-K}^{1/2} & 0_{(r-K) \times (T-r)} \\ 0_{(N-r) \times K} & 0_{(N-r) \times (r-K)} & 0_{(N-r) \times (T-r)} \end{bmatrix},$$

where $\Lambda_{r-K}^{1/2} = \operatorname{diag}(\sqrt{\lambda_{K+1}}, \ldots, \sqrt{\lambda_r})$ and the $\lambda_i$'s are the singular values of $Y'Y$. The columns of the $N \times N$ orthogonal matrix $W$ are the eigenvectors of $YY'$ (arranged according to the non increasing order of the associated eigenvalues) and the columns of the $T \times T$ orthogonal matrix $V$ are the eigenvectors of $Y'Y$. Note that some or all of the submatrices disappear when $r = K$, and/or $r = N$, and/or $r = T$.

Let the matrix $W$ be partitioned in such a way that it is conformable to the partitioned matrix $D$, that is

$$W = \begin{bmatrix} W(1:K,1:K) & W(1:K,K+1:r) & W(1:K,r+1:N) \\ W(K+1:r,1:K) & W(K+1:r,K+1:r) & W(K+1:r,r+1:N) \\ W(r+1:N,1:K) & W(r+1:N,K+1:r) & W(r+1:N,r+1:N) \end{bmatrix},$$

where $W(u:v,k:h)$ is used to denote the submatrix of $W$ which includes all entries $w_{ij}$ with $u \leq i \leq v$ and $k \leq j \leq h$. Multiplication by $D$ yields

$$WD = \begin{bmatrix} W(1:K,1:K)\Lambda_K^{1/2} & W(1:K,K+1:r)\Lambda_{r-K}^{1/2} & 0_{K\times(T-r)} \\ W(K+1:r,1:K)\Lambda_K^{1/2} & W(K+1:r,K+1:r)\Lambda_{r-K}^{1/2} & 0_{(r-K)\times(T-r)} \\ W(r+1:N,1:K)\Lambda_K^{1/2} & W(r+1:N,K+1:r)\Lambda_{r-K}^{1/2} & 0_{(N-r)\times(T-r)} \end{bmatrix}.$$

On the other hand, we have obviously $\hat{X}' = [u_1,\ldots,u_K]$ and multiplication of $V'$ by $\hat{X}'$ yields

$$V'\hat{X}' = [u_1,\ldots,u_T]'[u_1,\ldots,u_K] = \begin{bmatrix} I_K \\ 0_{(T-K)\times K} \end{bmatrix} = \begin{bmatrix} I_K \\ 0_{(r-K)\times K} \\ 0_{(T-r)\times K} \end{bmatrix}.$$

Some or all of the zero submatrices disappear when $r = K$ and/or $r = T$. Summing up, we obtain

$$Y\hat{X}' = WDV'\hat{X}' =$$

$$\begin{bmatrix} W(1:K,1:K)\Lambda_K^{1/2} & W(1:K,K+1:r)\Lambda_{r-K}^{1/2} & 0_{K\times(T-r)} \\ W(K+1:r,1:K)\Lambda_K^{1/2} & W(K+1:r,K+1:r)\Lambda_{r-K}^{1/2} & 0_{(r-K)\times(T-r)} \\ W(r+1:N,1:K)\Lambda_K^{1/2} & W(r+1:N,K+1:r)\Lambda_{r-K}^{1/2} & 0_{(N-r)\times(T-r)} \end{bmatrix} \begin{bmatrix} I_K \\ 0_{(r-K)\times K} \\ 0_{(T-r)\times K} \end{bmatrix}$$

$$= \begin{bmatrix} W(1:K,1:K)\Lambda_K^{1/2} \\ W(K+1:r,1:K)\Lambda_K^{1/2} \\ W(r+1:N,1:K)\Lambda_K^{1/2} \end{bmatrix} = \begin{bmatrix} W(1:K,1:K) \\ W(K+1:r,1:K) \\ W(r+1:N,1:K) \end{bmatrix} \Lambda_K^{1/2}.$$

It is easily recognized that

$$\begin{bmatrix} W(1:K,1:K) \\ W(K+1:r,1:K) \\ W(r+1:N,1:K) \end{bmatrix} = W_K,$$

so that the desired result

$$\hat{A} = W_K \Lambda_K^{1/2}$$

follows. □

### A.9. Proof of Theorem 6.1

*Proof.* The result for the mean is obvious. For the variance we have, for example,

$$\text{var}\{T(\hat{\gamma}_{rs}(0) - \tilde{\gamma}_{rs}(0))\} = \mathbb{E}\{(z_{rt_0} - \bar{z}_r)\omega_s + \omega_r(z_{st_0} - \bar{z}_s)\}^2$$

$$= \omega_s^2 \mathbb{E}(z_{rt_0} - \bar{z}_r)^2 + 2\omega_s\omega_r \mathbb{E}\{(z_{rt_0} - \bar{z}_r)(z_{st_0} - \bar{z}_s)\} + \omega_r^2 \mathbb{E}(z_{st_0} - \bar{z}_s)^2.$$

Then, by recalling that

$$\mathbb{E}(\tilde{\gamma}_{rs}(0)) = \gamma_{rs}(0) + O(T^{-1}),$$

(e.g. 17, vol. 2, p. 693, Formula 9.5.4) the result follows. The proof for $\text{var}\{T(\hat{\gamma}_{rs}(h) - \tilde{\gamma}_{rs}(h))\}$ is similar. □

*A.10. Proof of Theorem 6.2*

*Proof.* We have
$$\mathrm{Re}\{I_{rs}(\lambda_j) - \tilde{I}_{rs}(\lambda_j)\} = \frac{\omega_r \omega_s}{2\pi T} + A_{rs}(\lambda_j)$$

and
$$\mathrm{Im}\{I_{rs}(\lambda_j) - \tilde{I}_{rs}(\lambda_j)\} = B_{rs}(\lambda_j),$$

and obviously $\mathbb{E}A_{rs}(\lambda) = \mathbb{E}B_{rs}(\lambda) = 0$ for each $r$, $s$, and $\lambda$.

We evaluate now the second moments
$$\mathrm{var}\sqrt{T}A_{rs}(\lambda_j) = \frac{1}{(2\pi)^2 T}\mathbb{E}\{\sum_{t=1}^{T}\sum_{u=1}^{T}\cos\lambda_j(t_0 - t)\cos\lambda_j(t_0 - u)$$
$$\times (\omega_r(z_{st} - \bar{z}) + (z_{rt} - \bar{z})\omega_s)(\omega_r(z_{su} - \bar{z}) + (z_{ru} - \bar{z})\omega_s)\}$$
$$= \frac{1}{(2\pi)^2 T}\sum_{t=1}^{T}\sum_{u=1}^{T}\cos\lambda_j(t_0 - t)\cos\lambda_j(t_0 - u)$$
$$\times \{\omega_r^2 \gamma_{ss}(t-u) + \omega_s^2 \gamma_{rr}(t-u) + \omega_r\omega_s(\gamma_{rs}(t-u) + \gamma_{rs}(u-t))\}$$
$$= \frac{1}{(2\pi)^2 T}\sum_{h=-T+1}^{T-1}\{\omega_r^2\gamma_{ss}(h) + \omega_s^2\gamma_{rr}(h) + \omega_r\omega_s(\gamma_{rs}(h) + \gamma_{rs}(-h))\}\sum_v \cos\lambda_j v \cos\lambda_j(v+h),$$

where we have put $v = t_0 - t$ and are neglecting end effects. If $\lambda_j$ is a Fourier frequency,
$$\sum_v \cos\lambda_j v \cos\lambda_j(v+h) = \sum_v \cos\lambda_j v(\cos\lambda_j v\cos\lambda_j h - \sin\lambda_j v\sin\lambda_j h)$$
$$= \cos\lambda_j h\sum_v \cos^2\lambda_j v - \sin\lambda_j h\sum_v \cos\lambda_j v\sin\lambda_j v$$
$$= \frac{T}{2}\cos\lambda_j h,$$

and finally
$$\mathrm{var}\sqrt{T}A_{rs}(\lambda_j) = \frac{1}{(2\pi)^2 T}\sum_{h=-T+1}^{T-1}\{\omega_r^2\gamma_{ss}(h) + \omega_s^2\gamma_{rr}(h) + \omega_r\omega_s(\gamma_{rs}(h) + \gamma_{rs}(-h))\}\frac{T}{2}\cos\lambda_j h$$
$$\to \frac{1}{4\pi}\{f_{ss}(\lambda_j)\omega_r^2 + f_{rr}(\lambda_j)\omega_s^2 + \omega_r\omega_s(f_{rs}(\lambda_j) + f_{rs}(-\lambda_j))\}.$$

On repeating the same argument for $B_{rs}(\lambda_j)$ we obtain
$$\mathrm{var}\sqrt{T}B_{rs}(\lambda_j) = \frac{1}{(2\pi)^2 T}\sum_{t=1}^{T}\sum_{u=1}^{T}\sin\lambda_j(t_0-t)\sin\lambda_j(t_0-u)$$
$$\times\{\omega_r^2\gamma_{ss}(t-u) + \omega_s^2\gamma_{rr}(t-u) - \omega_r\omega_s(\gamma_{rs}(t-u) + \gamma_{rs}(u-t))\},$$

and since

$$\sum_v \sin\lambda_j v \sin\lambda_j(v+h) = \sum_v \sin\lambda_j v(\sin\lambda_j v\cos\lambda_j h + \cos\lambda_j v\sin\lambda_j h)$$

$$= \cos\lambda_j h \sum_v \sin^2\lambda_j v + \sin\lambda_j h \sum_v \sin\lambda_j v\cos\lambda_j v$$

$$= \frac{T}{2}\cos\lambda_j h,$$

it follows

$$\operatorname{var}\sqrt{T}B_{rs}(\lambda_j) \to \frac{1}{4\pi}\{f_{ss}(\lambda_j)\omega_r^2 + f_{rr}(\lambda_j)\omega_s^2 - \omega_r\omega_s(f_{rs}(\lambda_j) + f_{rs}(-\lambda_j))\}.$$

Noting that

$$\sqrt{T}\operatorname{Re}\{I_{rs}(\lambda_j(T)) - \tilde{I}_{rs}(\lambda_j(T))\} = \sqrt{T}A_{rs}(\lambda_j(T))$$

and

$$\sqrt{T}\operatorname{Im}\{I_{rs}(\lambda_j(T)) - \tilde{I}_{rs}(\lambda_j(T))\} = \sqrt{T}B_{rs}(\lambda_j(T))$$

the result follows. □

### A.11. Proof of Theorem 6.3

*Proof.* Since

$$\hat{F}_{rs}(\lambda) - \tilde{F}_{rs}(\lambda) = \frac{2\pi}{T}\sum_{j=-T/2}^{T/2} w_T(\lambda - \lambda_j)\{A(\lambda_j) + iB(\lambda_j)\} + \frac{\omega_r\omega_s}{2\pi T},$$

we have

$$\operatorname{var}\operatorname{Re}\{\hat{F}_{rs}(\lambda) - \tilde{F}_{rs}(\lambda)\} = (\frac{2\pi}{T})^2 \sum_j\sum_k w_T(\lambda-\lambda_j)w_T(\lambda-\lambda_k)\operatorname{cov}\{A_{rs}(\lambda_j)A_{rs}(\lambda_k)\}$$

and

$$\operatorname{var}\operatorname{Im}\{\hat{F}_{rs}(\lambda) - \tilde{F}_{rs}(\lambda)\} = (\frac{2\pi}{T})^2 \sum_j\sum_k w_T(\lambda-\lambda_j)w_T(\lambda-\lambda_k)\operatorname{cov}\{B_{rs}(\lambda_j)B_{rs}(\lambda_k)\},$$

and it may be easily seen that $\operatorname{cov}\{A_{rs}(\lambda_j)A_{rs}(\lambda_k)\}$ and $\operatorname{cov}\{B_{rs}(\lambda_j)B_{rs}(\lambda_k)\}$ tend to zero as $T \to \infty$ if $\lambda_j \neq \lambda_k$ while the results for $\lambda_j = \lambda_k$ are given in the previous theorem. As $T \to \infty$, the window bandwidth tends to zero, $(2\pi/T)\sum w_T(\lambda-\lambda_j)^2$ tends to $c_0 M$ and therefore $\operatorname{var}\{\operatorname{Re}(F_{rs}(\lambda) - \tilde{F}_{rs}(\lambda))\}$ has the same asymptotic behavior as $\frac{2\pi M}{T^2}\operatorname{var}(\sqrt{T}A_{rs}(\lambda))$ while $\operatorname{var}\{\operatorname{Im}(F_{rs}(\lambda) - \tilde{F}_{rs}(\lambda))\}$ has the same asymptotic behavior as $\frac{2\pi M}{T^2}\operatorname{var}(\sqrt{T}B_{rs}(\lambda))$. On substituting the expressions for the variance the thesis follows. □

### Acknowledgements

We thank an anonymous referee for careful comments and helpful suggestions.

*R.Baragona and F. Battaglia/Outliers in dynamic factor models* 431## References